\documentclass[11pt,leqno]{article}

\usepackage{amsthm, amsfonts, amssymb, epsfig, graphics,
amsmath, latexsym, eufrak}
\usepackage[latin1]{inputenc}\relax

\usepackage{aeguill}

\newcommand\uhu{\underline{\hat{u}}}

\newcommand{\CC}{{\mathbb C}}

\newcommand{\EE}{{\mathbb E}}
\newcommand{\FF}{{\mathbb F}}
\newcommand{\GG}{{\mathbb G}}

\newcommand{\MM}{{\mathbb M}}
\newcommand{\NN}{{\mathbb N}}

\newcommand{\PP}{{\mathbb P}}

\newcommand{\RR}{{\mathbb R}}

\newcommand{\VV}{{\mathbb V}}

\newcommand\cF{\mathcal{F}}

\newcommand\cH{\mathcal{H}}

\newcommand\cO{\mathcal{O}}

\def\UU {{\mathbb{U}}}

\def\bP{{\bf P}}
\def\ubP{\underline{\bf P}}

\def\bPi{\boldsymbol{\Pi}}

\newcommand\D{\partial }

\newcommand\Ker{\mathrm{Ker }}

\newcommand\uu{\underline{u }}

\newcommand\uR{\underline{R }}
\newcommand\uU{\underline{U }}
\newcommand\uV{\underline{V }}

\newcommand\hu{\hat{u}}
\newcommand\he{\hat{e}}
\newcommand\hw{\hat{w}}
\newcommand\hf{\hat{f}}
\newcommand\hv{\hat{v}}
\newcommand\hr{\hat{r}}

\newcommand\hV{\hat{V}}

\newcommand\tR{\tilde{R}}

\newcommand\tg{\tilde{g}}
\newcommand\tA{\tilde{A}}

\newcommand\tV{\widetilde{V }}

\newcommand\tu{\tilde{u }}

\newcommand\ur{\underline{r}}

\def\eps {\varepsilon}

\newtheorem{theo}{Theorem} [section]
\newtheorem{lem}[theo]{Lemma}
\newtheorem{cor}[theo]{Corollary}
\newtheorem{prop}[theo]{Proposition}

\newtheorem{hyp}[theo]{Assumption}
\newtheorem{rem}[theo]{Remark}

\numberwithin{equation}{section}

\newenvironment{proof-sketch}{\noindent{\bf Sketch of Proof}\hspace*{1em}}{\qed\bigskip\\}
\newenvironment{preuve} {\noindent {\sl{Proof. }}}
{\hfill $ \Box$ }

\begin{document}
\title
{ \LARGE \bf \sc Penalization approach for mixed hyperbolic systems
with constant coefficients satisfying a Uniform Lopatinski
Condition.}

\bigskip
\bigskip

 \author {{\sc B. Fornet}
 \footnote{LATP, Universit\'e de Provence,\footnotesize
39 rue Joliot-Curie, 13453 Marseille cedex 13, France.}
 \footnote{LMRS, Universit\'e de Rouen,\footnotesize
Avenue de l'Universit\'e, BP.12, 76801 Saint Etienne du Rouvray,
France.}
 }

\bigskip
\bigskip

\maketitle


\begin{abstract}
In this paper, we describe a new, systematic and explicit way of
approximating solutions of mixed hyperbolic systems with constant
coefficients satisfying a Uniform Lopatinski Condition via different
Penalization approaches.
\end{abstract}
\newpage

\section{Introduction.}
In this paper, we describe a new, systematic and explicit way of
approximating solutions of mixed hyperbolic systems with constant
coefficients satisfying a Uniform Lopatinski Condition via different
Penalization approaches. In applied Mathematics like, for instance,
in the study of fluids dynamics,  the method of penalization is used
to treat boundary conditions in the case of complex geometries. By
replacing the boundary condition by a singular perturbation of the
PDE extended to a larger domain, this method allows the construction
of an approximate, often more easily computable, solution.  We
consider mixed boundary value problems for hyperbolic systems:
$$\D_t+\sum_{j=1}^d A_j\D_j,$$ on $\{x_d\geq 0\},$ with boundary
conditions on $\{x_d=0\}.$  The $n\times n$ real valued matrices
$A_j$ are assumed constant. Of course, we assume the coefficients to
be constant as a first approach, aiming to generalize the results
obtained here in future works. We assume that the boundary
$\{x_d=0\}$ is noncharacteristic, which means that $det A_d\neq 0.$
We denote by
\newline $y:=(x_1,\ldots,x_{d-1})$ and $x:=x_d.$ The problem writes:
\begin{equation}
\left\{\begin{aligned}{} \label{deb}
& \cH u=f, \quad \{x>0\},\\
& \Gamma u|_{x=0}=\Gamma g ,\\
& u|_{t<0}=0 \quad ,
\end{aligned}\right.
\end{equation}
where the unknown  $u(t,x)\in\RR^n,$ $\Gamma:\RR^n\rightarrow \RR^p$
is linear and such that $rg \, \, \Gamma=p;$ which implies that
$\Gamma$ can be viewed as a $p\times n$ real valued constant matrix.
Let us fix $T>0$ once and for all for this paper. Let $\Omega_T^+$
denotes the set $[0,T]\times \RR^d_+$ and $\Upsilon_T$ denote the
set $[0,T]\times \RR^{d-1}.$ $f$ is a function in $H^k(\Omega_T^+),$
$g$ is a function in $H^k(\Upsilon_T)$,  where $k\geq 3$ or
$k=\infty,$ such that: $f|_{t<0}=0$ and $g|_{t<0}=0.$ We make
moreover the following Hyperbolicity assumption on $\cH:$
\begin{hyp}\label{vie}
For all $(\eta,\xi)\in \RR^{d-1}\times\RR-\{0\},$ the eigenvalues of
$$\sum_{j=1}^{d-1}\eta_jA_j+\xi A_d $$ are real, semi-simple and of
constant multiplicity.
\end{hyp}
Let us introduce now the frequency variable
$\zeta:=(\gamma,\tau,\eta),$ where $i\tau+\gamma,$ with $\gamma\geq
0,$ and $\tau\in\RR$ stands for the frequency variable dual to $t$
and $\eta=(\eta_1,\ldots,\eta_{d-1})$ where $\eta_j\in\RR$ is the
frequency variable dual to $x_j.$ We note:
$$A(\zeta):=-\left(A_d \right)^{-1}\left((i\tau+\gamma) Id+
\sum_{j=1}^{d-1} i\eta_j A_j\right).$$ Denote by $M$ a $N\times N,$
complex valued, matrix; $\EE_-(M)$[resp $\EE_+(M)$] is the linear
subspace generated by the generalized eigenvectors associated to the
eigenvalues of $M$ with negative [resp positive] real part. If $\FF$
and $\GG$ denote two linear subspaces of $\CC^N$ such that $\dim \FF
+ \dim \GG=N,$ $det(\FF,\GG)$  denotes the determinant obtained by
taking orthonormal bases in each space. Up to the sign, the result
is independent of the choice of the bases. We shall now explicit the
Uniform Lopatinski Condition assumption:
\begin{hyp}\label{vie2}
$(\cH,\Gamma)$ satisfies the Uniform Lopatinski Condition i.e for
all $\zeta$ such that $\gamma>0,$ there holds:
\begin{equation}\label{lopp}|det(\EE_{-}(A),\ker \Gamma)|\geq C>0.\end{equation}
\end{hyp}
The mixed hyperbolic system \eqref{deb} has a unique solution in
\newline $H^k(\Omega_T^+),$ and, since $\cH$ is hyperbolic with
constant multiplicity, for all $\gamma$ positive, the eigenvalues of
$A$ stay away from the imaginary axis. Moreover, as emphasized for
instance by Chazarain and Piriou in \cite{CP} and Métivier in
\cite{M}, there is a continuous extension of the linear subspace
$\EE_{-}(A)$ to $\{\gamma=0, (\tau,\eta)\neq 0_{\RR^d}\}$ that we
will denote by $\tilde{\EE}_{-}(A).$ $\tilde{\EE}_{+}(A)$ extends as
well continuously to $\{\gamma=0, (\tau,\eta)\neq 0_{\RR^d}\}$ and
we will denote $\tilde{\EE}_{+}(A)$ this extension. Moreover, there
holds:
$$\tilde{\EE}_{-}(A)\bigoplus \tilde{\EE}_{+}(A)=\CC^N.$$ We can refer the reader to \cite{CP}, \cite{GMWZ}, \cite{K},  or \cite{M} for detailed estimates concerning mixed hyperbolic problems satisfying a Uniform Lopatinski Condition.
Moreover, we can refer to \cite{MZ2} for the proof of the continuous
extension of the linear subspaces mentioned above in the
hyperbolic-parabolic framework.
\begin{rem}
As a consequence of the uniform Lopatinski condition, there holds, for all $\zeta\neq 0:$
$$rg \, \, \Gamma=p=\dim \tilde{\EE}_{-}(A(\zeta)).$$
\end{rem}

\subsection{A Kreiss Symmetrizer Approach.}

We will now describe a penalization method involving a Kreiss
Symmetrizer and a matrix constructed by Rauch in \cite{R}, in the
construction of our singular perturbation. Note well that we have
some freedom in both the choice of the Kreiss Symmetrizer and of
Rauch's matrix. Let us denote respectively by $\hu,$ $\hf,$ and
$\hat{g}$ the tangential Fourier-Laplace transform of $u,$ $f,$ and
$g.$ Since the Uniform Lopatinski Condition is holding for the mixed
hyperbolic system \eqref{deb}, there is, see \cite{MZ} a Kreiss
symmetrizer $S$ for the problem:
\begin{equation}\label{pb}
\left\{\begin{aligned}{}
& \D_x \hu= A \hu + \hf, \quad \{x>0\},\\
& \Gamma \hu|_{x=0}=\Gamma \hat{g} ,\\
\end{aligned}\right.
\end{equation}
That is to say there exists a matrix $S(\zeta),$ homogeneous of
order zero in $\zeta,$ $C^\infty$ in $\RR^+\times
\RR^d-\{0_{\RR^{d+1}}\}$ and there are $\lambda>0,$ $\delta>0$ and
$C_1$ such that:
\begin{itemize}
\item $S$ is hermitian symmetric.
\item $\Re \left(S A\right)\geq \lambda Id.$
\item $S\geq \delta Id-C_1\Gamma^*\Gamma.$
\end{itemize}
An algebraic result proved by Rauch in \cite{R} can be reformulated
as follow, and a proof is recalled in section \ref{prau}:
\begin{lem}\label{rau} There is a hermitian symmetric, uniformly definite positive, $N\times N$ matrix $B$
such that:
$$\ker \Gamma=\EE_+((S)^{-1}B).$$
Moreover $B$ depends smoothly of $\zeta.$
\end{lem}
\begin{rem}
This result is proved by constructing explicit matrices satisfying
the desired properties. Thus, it is not merely an existence result
and we can use the explicitly known matrix $B$ in our construction
of a penalization operator.
\end{rem}
\noindent Let us denote by $R:=B^{\frac{1}{2}}$ and $S_R:=R^{-1}S
R^{-1}.$ We will denote by $\PP^-$ the projector on $\EE_-(S_R)$
parallel to $\EE_+(S_R)$ and by $\PP^+$ the projector on
$\EE_+(S_R)$ parallel to $\EE_-(S_R);$ $\underline{\PP}^-$ and
$\underline{\PP}^+$ denoting the associated Fourier multiplier. We
recall that, denoting by $\cF$ the tangential Fourier transform, the
Fourier multiplier $\underline{\PP}^-(\D_t,\D_y,\gamma)$ [resp
$\underline{\PP}^+(\D_t,\D_y,\gamma)$] is then defined, for all
$w\in H^k(\RR^{d+1}),$ and $\gamma>0,$ by:
$$\cF\left(\underline{\PP}^-(\D_t,\D_y,\gamma) w\right)=\PP^-(\zeta)\cF(w),$$ [resp $$\cF\left(\underline{\PP}^+(\D_t,\D_y,\gamma) w\right)=\PP^+(\zeta)\cF(w)],$$ in the future we will rather write:
$$\cF\left(\underline{\PP}^\pm(\D_t,\D_y,\gamma) w\right)=\PP^\pm(\zeta)\cF(w).$$
We fix, once and for all, $\gamma>0$ big enough. Let us consider then the solution $\underline{u}^\eps$ of the well-posed Cauchy problem on the whole space \eqref{penak}:
\begin{equation}\label{penak}
\left\{\begin{aligned}{}
&\cH \underline{u}^\eps+ \frac{1}{\eps} \MM \underline{u}^\eps \textbf{1}_{x<0} =f\textbf{1}_{x>0}+ \frac{1}{\eps}\theta \textbf{1}_{x<0}, \quad \{x\in\RR\},\\
&\underline{u}^\eps|_{t<0}=0,\\
\end{aligned}\right.
\end{equation}
where $$\MM:= -e^{\gamma t} A_d \underline{S}^{-1} \underline{R} \underline{\PP}^-
\underline{R} e^{-\gamma t},$$ $$\theta:=-e^{\gamma t} A_d
\underline{S}^{-1}\underline{R}\underline{\PP}^- \underline{\Gamma}
\tg,$$ and
$\underline{S}(\D_t,\D_y)$ [resp $\underline{R}(\D_t,\D_y)$] denotes
the Fourier multiplier associated to $S(\zeta)$ [resp $R(\zeta)$].
Let us define $\tilde{g}$ by:
$$\tg:=e^{-x^2}g.$$
In what follows, $\hat{g}$ will denote the Fourier-Laplace transform
of $\tg.$ Let us denote by
$$\underline{\tu}:=\underline{u}^-\textbf{1}_{x<0}+u\textbf{1}_{x\geq 0}=\underline{u}^-\textbf{1}_{x\leq 0}+u\textbf{1}_{x>0}.$$
$u$ denotes the solution of \eqref{deb}, and thus belongs to
$H^{k}(\Omega_T^+).$ $\underline{u}^-$ is a function belonging to
$H^{k}(\Omega_T^-)$ and such that $\underline{u}^-|_{x=0}=u|_{x=0}.$
More precisely, $\underline{u}^-$ can be computed by: $e^{\gamma
t}\cF^{-1}\left(R^{-1}(\underline{\hv}^-+\PP^-\Gamma\hat{g})\right),$
where $\underline{\hv}^-$ is the solution of the problem:
\begin{equation*}
\left\{\begin{aligned}{} &S_R \D_x \underline{\hv}^- - \PP^+ S_R A_R
\underline{\hv}^-=\PP^+ S_R A_R \PP^-\Gamma \hat{g}, \quad \{x<0\},\\
&\underline{\hv}^-|_{x=0}=\PP^+R\hu|_{x=0},\\
\end{aligned}\right.
\end{equation*}
and $\hu$ denotes the Fourier-Laplace transform of the solution $u$
of \eqref{deb}.
\begin{theo}\label{kreiss}
For all $k\geq 3,$ if $f\in H^{k}(\Omega_T^+)$ and $g\in
H^{k}(\Upsilon_T),$ then there holds:
$$\|\uu^\eps-\underline{u}^-\|_{H^{k-3}(\Omega_T^-)}+ \|\uu^\eps-u\|_{H^{k-3}(\Omega_T^+)}=\cO(\eps),$$
where $\uu^\eps$ denotes the solution of the Cauchy problem
\eqref{penak} and $u$ denotes the solution of the mixed hyperbolic
problem \eqref{deb}. If $g=0$ then:
$$\|\uu^\eps-\underline{u}^-\|_{H^{k-\frac{3}{2}}(\Omega_T^-)}+ \|\uu^\eps-u\|_{H^{k-\frac{3}{2}}(\Omega_T^+)}=\cO(\eps).$$
\end{theo}
Of course, since $\uu^\eps$ is defined for all $\{x\in\RR\},$ its
limit as $\eps\rightarrow 0^+,$ $\underline{\tu}$ is can be viewed
as an "extension" of $u$ on the fictive domain $\{x<0\}.$ The
"extension" resulting from our method of penalization gives a
continuous $\underline{\tu}$ across $\{x=0\},$ while the method used
in \cite{BR} gave simply: $\underline{\tu}|_{x<0}=0.$ We have the
following Corollaries:
\begin{cor}
Assume for example that $f\in H^\infty(\Omega_T^+)$ and\newline
$g\in H^\infty(\Upsilon_T)$ then
$$\|\uu^\eps-u\|_{H^{s}(\Omega_T^+)}=\cO(\eps);\quad \forall s>0.$$
\end{cor}
\begin{cor}\label{cor}
If $f$ belongs to $L^2(\Omega_T^+)$ and $g=0$ then:
$$\lim_{\eps\rightarrow 0^+} \|\uu^\eps-\underline{\tu}\|_{L^2(\Omega_T)}=0.$$
\end{cor}
One of the interest of this first approach lies in the rate of
convergence of $\uu^\eps$ towards $u$. Indeed, in general, a
boundary layer will form near the boundary in this kind of singular
perturbation problem. For example in the paper  by Bardos and Rauch
\cite{BR}, as confirmed by Droniou \cite{D}, a boundary layer forms.
It is also the case in \cite{paccou}, as analyzed in our Appendix.
There are also boundary layers phenomena in the parabolic context:
see the approach proposed by Angot, Bruneau and Fabrie \cite{ABF}
for instance. However, surprisingly, and like in the penalization
method proposed by Fornet and Guès in \cite{FG1}, our method allows
the convergence to occur without formation of any boundary layer on
the boundary. As a result, this leads to the kind of sharp stability
estimate given in Theorem \ref{kreiss}. These results concern the
case where $f$ and $g$ are sufficiently regular. The reason is that
we construct an approximate solution. In the case of $g$ only in
$L^2(\Upsilon_T),$ such a simple treatment does not work. However,
let $\delta>0$ be given. If we approximate $f$ and $g$ by smooth
functions $f_\nu\in H^\infty(\Omega_T^+)$ and $g_\nu\in
H^\infty(\Upsilon_T)$ such that
$\|f-f_\nu\|_{L^2(\Omega_T^+)}<\delta$ and
$\|g-g_\nu\|_{L^2(\Upsilon_T)}<\delta,$ by the uniform Lopatinski
condition, we get: $$\|u_\nu-u\|_{L^2(\Omega_T^+)}<C\nu,$$ where
$u_\nu$ is the solution of the mixed hyperbolic problem \eqref{deb}
with data $f_\nu$ and $g_\nu.$ We can now apply Corollary \ref{cor}
to $u_\nu,$ and obtain by penalization a sequence $u_\nu^\eps$ in
$L^2(\Omega_T)$ such that: $\lim_{\eps\rightarrow 0^+}
u_\nu^\eps=u_\nu$ in $L^2(\Omega_T^+).$ Finally, by choosing, $\eps$
sufficiently small, we get
$\|u-u_\nu^\eps\|_{L^2(\Omega_T^+)}<2C\delta.$ By choosing $\eps$
and $\nu$ as functions of $\delta,$ and noting
$u_{(\delta)}=u_{\nu(\delta)}^{\eps(\delta)},$ we have:
\begin{equation}\label{delta}
\begin{aligned}{} & \lim_{\delta\rightarrow 0^+}
\|u_{(\delta)}-u\|_{L^2(\Omega_T^+)}=0.\\
\end{aligned}
\end{equation}

\subsection{A second Approach.}

In the first approach we have just introduced, it is necessary to
compute a Kreiss's Symmetrizer and a Rauch's matrix. In view of
future numerical applications, we will now introduce another method
preventing the computation of these matrices. The price to pay is
that we need the preliminary computation of $v,$ which is by
definition the solution of the Cauchy problem on the free space:
\begin{equation}
\left\{\begin{aligned}{} \label{cauch}
& \cH v=f, \quad (t,y,x)\in \Omega_T,\\
& v|_{t<0}=0 \quad \forall(y,x)\in\RR^d.
\end{aligned}\right.
\end{equation}

Let us denote $\bP^-(\zeta)$ the spectral projector on
$\tilde{\EE}_-(A(\zeta))$ parallel to $\tilde{\EE}_+(A(\zeta)),$ and
$\bP^+(\zeta)$ the spectral projector on $\tilde{\EE}_+(A(\zeta))$
parallel to $\tilde{\EE}_-(A(\zeta)).$ Let us introduce
$\ubP^\pm(\D_t,\D_y,\gamma),$ the Fourier multiplier associated to
$\bP^\pm(\zeta).$ Let us denote by $\bPi$ the projector on
$\tilde{\EE}_-(A(\zeta))$ parallel to $\Ker \Gamma,$ which has a
sense because of the Uniform Lopatinski Condition and denote
$\underline{\bPi}$ the associated Fourier multiplier. We
define then $\tilde{h}$ by:
$$\tilde{h}:=e^{-x^2}\left(\ubP^-(e^{-\gamma t}v|_{x=0})+\underline{\bPi} e^{-\gamma t} (g-v|_{x=0})\right),$$
where $g$ denotes the function involved in the boundary condition of
the mixed hyperbolic problem \eqref{deb}. Now, let us consider the
following singularly perturbed Cauchy problem on the whole space:
\begin{equation}
\left\{\begin{aligned}{} \label{pena}
& \cH u^\eps+\frac{1}{\eps} A_d e^{\gamma t} \ubP^- e^{-\gamma t} u^\eps \textbf{1}_{x<0}=f\textbf{1}_{x>0}+\frac{1}{\eps} A_d e^{\gamma t} \tilde{h} \textbf{1}_{x<0},\\
& u^\eps|_{t<0}=0 \quad .
\end{aligned}\right.
\end{equation}
Let us denote by
$$\tu:=u^-\textbf{1}_{x<0}+u\textbf{1}_{x\geq
0}=u^-\textbf{1}_{x\leq 0}+u\textbf{1}_{x>0}.$$ $u$ denotes the
solution of \eqref{deb} thus belonging to $H^{k}(\Omega_T^+)$ and
$u^-$ is a function belonging to $H^{k}(\Omega_T^-)$ and such that
$u^-|_{x=0}=u|_{x=0}.$ More precisely, $u^-$ can be computed by:
$e^{\gamma t}\cF^{-1}(\cF(\tilde{h})+\hv^-),$ where $\hv^-$ is the
solution of the problem:
\begin{equation*}
\left\{\begin{aligned}{}
& \D_x (\bP^+\hv^-)-A (\bP^+\hv^-)= 0, \quad \{x<0\},\\
& \bP^+\hv^-|_{x=0}= \bP^+\hu|_{x=0}.\\
\end{aligned}\right.
\end{equation*}
and $\hu$ denotes the Fourier-Laplace transform of the solution $u$
of \eqref{deb}. The problem \eqref{pena} is well-posed and, for all
$\eps>0,$ there exists a unique \newline $u^\eps\in H^k(\Omega_T)$ solution.
We will fix $\gamma$ adequately big beforehand.  We observe then the
following result:
\begin{theo}\label{maintheo}
For all $k\geq 3,$ if $f\in H^k(\Omega_T^+)$ and $g\in
H^k(\Upsilon_T),$ then there holds:
$$\|u^\eps-u^-\|_{H^{k-3}(\Omega_T^-)}+ \|u^\eps-u\|_{H^{k-3}(\Omega_T^+)}=\cO(\eps),$$
where $u^\eps$ denotes the solution of the Cauchy problem
\eqref{pena} and $u$ denotes the solution of the mixed hyperbolic
problem \eqref{deb}.
\end{theo}
The singular perturbation involved in the definition of $u^\eps$
does not depend either of Kreiss's Symmetrizer or Rauch's  matrix.
As a result, for this method of penalization far less computations
are necessary in order to obtain our singular perturbation. Note
well that the proof of the energy estimates in Theorem
\ref{maintheo} is completely different from the proof of the energy
estimates in Theorem \ref{kreiss}. Indeed, for our first approach
our singularly perturbed problem was treated as a Cauchy problem,
contrary to our second approach where it was interpreted as a
transmission problem.
\begin{cor}
Assume for example that $f\in H^\infty(\Omega_T^+)$ and\newline
$g\in H^\infty(\Upsilon_T)$ then
$$\|u^\eps-u\|_{H^{s}(\Omega_T^+)}=\cO(\eps);\quad \forall s>0.$$
\end{cor}
Of course, we see that the same problem of regularity arises in
Theorem \ref{maintheo} and Theorem \ref{kreiss}. However, by a
simple density argument, we can also prove here the exact analogous
of \eqref{delta}.
\begin{rem}
In the case where $f=0,$ then the solution $v$ of \eqref{cauch} is $v=0$ and thus, the perturbed cauchy problem \eqref{pena} rewrites:
\begin{equation*}
\left\{\begin{aligned}{}
& \cH u^\eps+\frac{1}{\eps} A_d e^{\gamma t} \ubP^- e^{-\gamma t} u^\eps \textbf{1}_{x<0}=\frac{1}{\eps} A_d e^{\gamma t} e^{-x^2}\left(\underline{\bPi} e^{-\gamma t} g\right) \textbf{1}_{x<0}, \quad \{x\in\RR\},\\
& u^\eps|_{t<0}=0 \quad .
\end{aligned}\right.
\end{equation*}
\end{rem}

\section{Underlying approach leading to the proof of\newline Theorem \ref{kreiss}.}

\subsection{Some preliminaries.}
Since the Uniform Lopatinski Condition holds, there is $S,$
homogeneous of order zero in $\zeta,$ and such that there are
$\lambda>0,$ $\delta>0$ and $C_1$ and there holds:
\begin{itemize}
\item $S$ is hermitian symmetric.
\item $\Re \left(S A\right)\geq \lambda Id.$
\item $S\geq \delta Id-C_1\Gamma^*\Gamma.$
\end{itemize}
$S$ is then called a Kreiss Symmetrizer for the problem:
\begin{equation}\label{pbb}
\left\{\begin{aligned}{}
& \D_x \hu= A \hu + \hf, \quad \{x>0\},\\
& \Gamma \hu|_{x=0}=\Gamma \hat{g} ,\\
\end{aligned}\right.
\end{equation}
where $\hf$ and $\hat{g}$ denotes respectively the Fourier-Laplace
transforms of $f$ and $\tg;$ and $\hu$ denotes the Fourier-Laplace
transform of the solution $u$ of the well-posed mixed hyperbolic
problem \eqref{deb}. $\hu$ is also solution, for all fixed
$\zeta\neq 0$ of the following equation:
\begin{equation}
\left\{\begin{aligned}{} \label{debm}
& S \D_x \hu= S A \hu + S (A_d)^{-1}\hf, \quad \{x>0\},\\
& \Gamma \hu|_{x=0}=\Gamma \hat{g} ,\\
\end{aligned}\right.
\end{equation}

\begin{rem}
Following our current assumptions, $\Gamma$ is independent of
$\zeta\neq 0,$ however, more general boundary conditions, of the
form: $$\Gamma(\zeta) \hu|_{x=0}=\Gamma(\zeta) \hat{g},$$ can be
treated. It would imply taking as boundary condition for
\eqref{deb}:
$$\Gamma_\gamma u|_{x=0}=\Gamma_\gamma g,$$ with for $\gamma$ big
enough,
$$\Gamma_\gamma:=\underline{\Gamma}(\D_t,\D_y) e^{-\gamma t},
$$
where, $\underline{\Gamma}(\D_t,\D_y)$ denotes the Fourier
multiplier associated to $\Gamma(\zeta),$ that is to say is defined
by:
$$\cF(\underline{\Gamma}(\D_t,\D_y) u)=\Gamma(\zeta) \cF(u).$$
\end{rem}
Referring for example to \cite{CP} and \cite{K}, Kreiss has proved
that the existence of a Kreiss symmetrizer for the symbolic equation
is sufficient to prove the well-posedness of the associated
pseudodifferential equation (here \eqref{deb}). Indeed, multiplying
by $\hu$ and integrating by parts the equation:
$$S \D_x \hu= S A \hu + S (A_d)^{-1}\hf$$ leads
to the desired a priori estimates. For all $\zeta\neq0,$ $S(\zeta)$
is hermitian symmetric and definite positive on $\ker \Gamma.$ Let
us sum up the properties crucial in the proof of the well-posedness
of our problem:
\begin{prop}\label{sym}
For all $\zeta=(\tau,\gamma,\eta)$ such that
$\tau^2+\gamma^2+\sum_{j=1}^{d+1} \eta_j^2=1,$ there holds:
\begin{itemize}
\item $S(\zeta)$ is hermitian symmetric.
\item $\Re \left(S A\right)(\zeta):=\frac{1}{2}(S A+ (SA)^*)(\zeta)$ is positive definite.
\item $-S(\zeta)$ is definite negative on $\ker \Gamma$ and $\ker
\Gamma$ is of same dimension as the number of negative eigenvalues
in $-S(\zeta).$
\end{itemize}
\end{prop}
\noindent Note that, by homogeneity of $S,$ it is equivalent for the
properties in Proposition \ref{sym} to hold for $|\zeta|=1$ or for
$|\zeta|>0.$ As a consequence of the first point and third point of
Proposition \ref{sym}, the Lemma \ref{rau} applies and gives a
matrix $B$ such that: $\ker\Gamma =\EE_+(S^{-1}B).$ In the sequel,
such a matrix $B$ is fixed once for all.

\noindent The following chapter contains a proof of Lemma \ref{rau}
assorted of a detailed construction of $B.$

\subsection{Detailed proof of Lemma \ref{rau}: Construction of the matrices $B$ solving Lemma
\ref{rau}.}\label{prau}

As we will emphasize in next chapter, Lemma \ref{rau} is a crucial
feature in our first method of Penalization. The aim of this chapter
is to give a more complete proof rather than simply recalling
Rauch's result and, in the process, to precise how the matrices $B$
solving Lemma \ref{rau} are constructed. For all $\zeta\neq 0,$
$S(\zeta)$ is hermitian symmetric, uniformly definite positive on
$\tilde{\EE}_+(A(\zeta)),$ and uniformly definite negative on
$\tilde{\EE}_-(A(\zeta));$ as a consequence, $S(\zeta)$ keeps
exactly $p$ positive eigenvalues and $N-p$ negative eigenvalues for
all $\zeta\neq0.$ Basically, knowing that $S$ is uniformly definite
positive on $\ker \Gamma;$ we search to express $\ker\Gamma$ in a
way involving $S.$ Consider $q\in\ker\Gamma,$ since, for all
$\zeta\neq 0,$ $\EE_-(S(\zeta))\bigoplus\EE_+(S(\zeta))=\CC^N,$ we
can split $q$ in: $$q:=q^+ + q^-$$ with $q^+\in\EE_+(S(\zeta))$ and
$q^-\in\EE_-(S(\zeta)).$ \newline Since $\dim ker \Gamma=\dim
\EE_+(S(\zeta))=p,$ these two linear subspaces are in bijection. Let
us give the two main ideas behind this proof: one idea is to detail
the bijection between $q\in ker \Gamma$ and $q^+\in\EE_+(S(\zeta))$
as it satisfies some constraints, the other is to come down to the
model case where the eigenvalues of $S$ are either $1$ or $-1.$ Let
us denote:
\begin{equation*}
\begin{aligned}{}
& \tilde{S}^{-1} = \left[\begin{array}{cc}
-Id_{N-p} & 0\\
0 & Id_{p}\\
\end{array}\right],
\end{aligned}
\end{equation*}
In a first step, we will prove the following result:
\begin{prop}\label{Id}
If we assume that $\VV$ is a linear subspace of $\CC^N$ of dimension
$p,$ and that there is $C>0$ such that, for all $q\in\VV,$ there
holds:
$$\langle \tilde{S}^{-1}q,q\rangle\geq C\langle q, q \rangle,$$
then the two following equivalent properties hold:
\begin{itemize}
\item There is a hermitian symmetric, positive definite matrix  $\tilde{\uR},$ such that:
$$\left[q\in\VV\right] \Leftrightarrow \left[\tilde{\uR}^{-1}q\in \EE_+(\tilde{\uR} \tilde{S}\tilde{\uR})\right],$$
which is equivalent to: $$\VV=\EE_+(\tilde{\uR}^2 \tilde{S}).$$

\item There is a hermitian symmetric, positive matrix $\tilde{R},$
such that:
$$\left[q\in\VV\right] \Leftrightarrow \left[\tilde{R}q\in \EE_+(\tilde{R} \tilde{S}^{-1}\tilde{R})\right],$$
which is equivalent to:
$$\VV=\EE_+(\tilde{S}^{-1}\tilde{R}^2).$$
\end{itemize}
Moreover, we can link the two properties by taking:
$$\tilde{R}^2=\tilde{S}\underline{\tilde{R}}^2 \tilde{S}.$$
\end{prop}

\begin{preuve}
In this proof, we will show how to construct some matrices
$\tilde{R}$ satisfying the required properties. There is a
$(N-p)\times p$ matrix $\aleph$ of rank $N-p$ such that
$\|\aleph\|\leq 1$ and:
$$\VV = \{q\in\CC^N, \quad q^-=\aleph q^+\},$$ where $q^+$[resp
$q^-$] denotes the projector on $\EE_+((\tilde{S})^{-1})$ [resp
$\EE_-((\tilde{S})^{-1})$]parallel to $\EE_-((\tilde{S})^{-1})$
[resp $\EE_+((\tilde{S})^{-1})$]. Indeed, $\dim \VV=p=\dim
\EE_+((\tilde{S})^{-1}),$ and
$\CC^N=\EE_-((\tilde{S})^{-1})\bigoplus \EE_+((\tilde{S})^{-1}).$
Moreover, there ic $C>0$ such that, for all $q\in\VV,$ there holds:
$$\langle (\tilde{S})^{-1}q,q \rangle= -\langle q^-,q^- \rangle+ \langle q^+,q^+ \rangle \geq C\langle q, q \rangle.$$
and thus $$|q^+|^2-|\aleph q^+|^2\geq C|q|^2,$$ which implies that
$\|\aleph\|<1.$ We will show now that, for $\tR$ constructed as
follow:
\begin{equation*}
\begin{aligned}{}
& \tR = \left[\begin{array}{cc}
Id_{N-p} & -\aleph\\
-\aleph^* & Id_{p}\\
\end{array}\right],
\end{aligned}
\end{equation*}
there holds:
$$\left[q\in\VV\right] \Leftrightarrow \left[\tilde{R}q\in \EE_+(\tilde{R} \tilde{S}^{-1}\tilde{R})\right].$$ First,we see that the constructed
$\tR$ is trivially hermitian symmetric and positive definite since
$\|\aleph\|< 1.$ First, we have:
\begin{equation*}
\begin{aligned}{}
& \tR \tilde{S}^{-1} \tR= \left[\begin{array}{cc}
-Id_{N-p}+NN^* & 0\\
0 & Id_{p}-N^*N\\
\end{array}\right],
\end{aligned}
\end{equation*}
and
\begin{equation*}
\begin{aligned}{}
& \tR q = \left(\begin{array}{c}
q^- -\aleph q^+\\
-\aleph^* q^-+ q^+\\
\end{array}\right).
\end{aligned}
\end{equation*}
Thus, since $\|\aleph\|<1,$ there holds: $$\left[\tilde{R}q\in
\EE_+(\tilde{R} \tilde{S}^{-1}\tilde{R})\right] \Leftrightarrow
\left[ q^- -\aleph q^+=0 \right]\Leftrightarrow
\left[q\in\VV\right].$$ We will now prove that we have:
$$(\tilde{R})^{-1}\EE_+(\tilde{R} \tilde{S}^{-1}\tilde{R})=\EE_+(\tilde{S}^{-1}\tilde{R}^2).$$
Since $\tilde{R} \tilde{S}^{-1}\tilde{R}$ is hermitian symmetric,
the linear subspace $\EE_+(\tilde{R} \tilde{S}^{-1}\tilde{R})$ is
generated by the eigenvectors of $\tilde{R} \tilde{S}^{-1}\tilde{R}$
associated to positive eigenvalues. A basis of
$(\tilde{R})^{-1}\EE_+(\tilde{R} \tilde{S}^{-1}\tilde{R})$ is thus
given by $((\tilde{R})^{-1}v_j)_j$ where $v_j$ denotes an
eigenvector of $\tilde{R} \tilde{S}^{-1}\tilde{R}$ associated to a
positive eigenvalue $\lambda_j$. We have:
$$\tilde{R} \tilde{S}^{-1}\tilde{R}v_j=\lambda_j v_j.$$ Let us
denote $w_j=(\tilde{R})^{-1}v_j,$ we have then:
$$\tilde{R} \tilde{S}^{-1}\tilde{R}^2 w_j=\lambda_j \tilde{R} w_j\Leftrightarrow \tilde{S}^{-1}\tilde{R}^2 w_j=\lambda_j w_j.$$
As a result, $w_j$ is an eigenvector of ${S}^{-1}\tilde{R}^2$
associated to the eigenvalue $\lambda_j$ hence we obtain that:
$$(\tilde{R})^{-1}\EE_+(\tilde{R} \tilde{S}^{-1}\tilde{R})=\EE_+(\tilde{S}^{-1}\tilde{R}^2).$$
We can also prove, the same way, that:
$$\tilde{\uR}\EE_+(\tilde{\uR} \tilde{S}\tilde{\uR})=\EE_+(\tilde{\uR}^2 \tilde{S}).$$
Now, taking $$\tilde{R}^2=\tilde{S}\underline{\tilde{R}}^2
\tilde{S},$$ we can check that:
$$\EE_+(\tilde{S}^{-1}\tilde{R}^2)=\EE_+(\tilde{\uR}^2 \tilde{S}),$$
which concludes the proof.
\end{preuve}

Lemma \eqref{r} is a Corollary of the following Proposition:
\begin{prop}\label{prop2}
If $S^{-1}$ denotes a smooth in $\zeta\neq0,$ matrix-valued function
in the space of hermitian symmetric matrices with $p$ positive
eigenvalues and $N-p$ negative eigenvalues and $\ker \Gamma$ denotes
a linear subspace of dimension $p$ and there is $C>0$ such that, for
all $q\in\ker\Gamma,$ there holds:
$$\langle S^{-1}q,q\rangle\geq C\langle q, q \rangle,$$
then the two following equivalent properties hold:
\begin{itemize}
\item There is a smooth in $\zeta\neq0,$ matrix-valued function
$\uR,$ in the space of hermitian symmetric, positive matrices such
that:
$$\left[q\in\Ker \Gamma\right] \Leftrightarrow \left[\forall \zeta\neq0, \quad \uR^{-1}(\zeta)q\in \EE_+(\uR(\zeta) S(\zeta)\uR(\zeta))\right],$$
which is equivalent to: $$\forall \zeta\neq 0, \quad \Ker
\Gamma=\EE_+(\uR^2(\zeta) S(\zeta)).$$
\item There is a smooth in $\zeta\neq0,$ matrix-valued function $R,$
in the space of hermitian symmetric, positive matrices such that:
$$\left[q\in\Ker \Gamma\right] \Leftrightarrow \left[\forall \zeta\neq0, \quad R(\zeta)q\in \EE_+(R(\zeta) S^{-1}(\zeta)R(\zeta))\right],$$
which is equivalent to: $$\forall \zeta\neq0, \quad \Ker
\Gamma=\EE_+(S^{-1}(\zeta)R^2(\zeta)).$$
\end{itemize}
Moreover, for all $\zeta\neq 0,$ these two properties can be linked
by taking:
$$(R(\zeta))^2=S(\zeta)(\uR(\zeta))^2 S(\zeta).$$
\end{prop}

\begin{preuve}
We will show here that Proposition \ref{prop2} can be deduced from
Proposition \ref{Id}. For all $\zeta\neq 0,$ $S(\zeta)$ is a
hermitian symmetric matrix, moreover $S$ depends smoothly of
$\zeta.$ As a consequence $S^{-1}$ is also  a hermitian symmetric
matrix depending smoothly of $\zeta,$ and as such,
 there is a nonsingular matrix $V$ such that:
$$\tilde{S}^{-1}=V^* \left(S^{-1}\right) V.$$ Let us denote $\Lambda$ the diagonalized  version of $S^{-1}$
with eigenvalues sorted by increasing order, then there is $Z$
depending smoothly of $\zeta$ such that, for all $\zeta\neq 0,$ we
have:
$$Z^*(\zeta)=Z^{-1}(\zeta),$$ and
$$\Lambda(\zeta)=Z^*(\zeta) \left(-S^{-1}\right)(\zeta)
Z(\zeta).$$ As a consequence, $V$ depends smoothly of $\zeta$ since,
for all $\zeta\neq 0$:
$$V(\zeta)= (\underline{\Lambda}(\zeta))^{-\frac{1}{2}} Z(\zeta),$$ where $\underline{\Lambda}$ is the diagonal matrix obtained by taking
the absolute value of each eigenvalue of $\Lambda.$ For the sake of
simplicity, let us omit the dependence in $\zeta.$ Now, for all
$q\in V^{-1}\ker\Gamma,$ there is $C>0,$ such that:
$$\langle \tilde{S}^{-1}q,q\rangle=\langle V^*S^{-1}Vq,q\rangle=\langle S^{-1}(Vq),(Vq)\rangle\geq C \langle (Vq),(Vq)\rangle. $$
Moreover $V$ is nonsingular, thus there is $C'>0$, such that, for
all \newline $q\in V^{-1}\ker\Gamma,$ there holds:
$$\langle \tilde{S}^{-1}q,q\rangle\geq C'\langle q,q\rangle.$$
Moreover $\dim V^{-1}\ker\Gamma=p,$ using Proposition \ref{Id}, for
all fixed $\zeta\neq 0,$ there is a hermitian symmetric, positive
definite matrix  $\tilde{\uR}(\zeta),$ such that:
$$V^{-1}(\zeta)\ker\Gamma=\EE_+((\tilde{\uR}(\zeta))^2
\tilde{S}(\zeta))=\tilde{\uR}(\zeta)\EE_+(\tilde{\uR}(\zeta)
\tilde{S}(\zeta)\tilde{\uR})(\zeta).$$ We will now prove that we can
construct $\tilde{\uR}$ depending smoothly of $\zeta.$  First there
is a $(N-p)\times p$ matrix $\aleph$ of rank $N-p$, depending
smoothly of $\zeta,$ such that fore all $\zeta\neq 0$
$\|\aleph(\zeta)\|\leq 1$ and:
$$V^{-1}(\zeta)\ker\Gamma= \{q\in\CC^N, \quad q^-=\aleph(\zeta) q^+\},$$ where $q^+$[resp
$q^-$] denotes the projector on $\EE_+((\tilde{S})^{-1})$ [resp
$\EE_-((\tilde{S})^{-1})$]parallel to $\EE_-((\tilde{S})^{-1})$
[resp $\EE_+((\tilde{S})^{-1})$]. $\underline{\tilde{R}}$ is given,
for all $\zeta\neq 0,$ by:
$$\underline{\tilde{R}}(\zeta)=\sqrt{\tilde{S}^{-1}(\zeta) \tilde{R}^2 (\zeta)\tilde{S}^{-1}(\zeta)},$$
with $\tR$ given by:
\begin{equation*}
\begin{aligned}{}
& \tR(\zeta) = \left[\begin{array}{cc}
Id_{N-p} & -\aleph(\zeta)\\
-\aleph^*(\zeta) & Id_{p}\\
\end{array}\right].
\end{aligned}
\end{equation*}
Since $\tilde{S}^{-1}=V^* \left(S^{-1}\right) V,$ there holds:
$\tilde{S}=V^* S V,$ and, as a consequence:
$$(V\tilde{\uR})^{-1}\ker\Gamma=\EE_+(\tilde{\uR}V^*S V\tilde{\uR}).$$
As $\tilde{\uR}V^*S V\tilde{\uR}$ is hermitian symmetric, a basis of
the linear subspace $\EE_+(\tilde{\uR}V^*S V\tilde{\uR})$ is given
by the eigenvectors of $\tilde{\uR}V^*S V\tilde{\uR}$ associated to
positive eigenvalues. This leads us to consider
$v_j=(V\tilde{\uR})^{-1} u_j$ satisfying:
$$\tilde{\uR}V^*S V\tilde{\uR}v_j=\lambda_jv_j.$$
We have:
$$\tilde{\uR}V^*S V\tilde{\uR}(V\tilde{\uR})^{-1} u_j=\lambda_j(V\tilde{\uR})^{-1} u_j.$$
hence:
$$(V\tilde{\uR})\tilde{\uR}V^*S  u_j=\lambda_j u_j.$$
Since
$(V\tilde{\uR})\tilde{\uR}V^*=(\tilde{\uR}V^*)^*(\tilde{\uR}V^*)$ is
hermitian symmetric and positive definite, we can then define its
square root. We define $\uR$ by:
$$\uR=\sqrt{(\tilde{\uR}V^*)^*(\tilde{\uR}V^*)}.$$
Since both $\tilde{\uR}$ and $V$ depends smoothly of $\zeta,$ so
does $\uR.$ Moreover, there holds:
$$\uR^2 S  u_j=\lambda_j u_j,$$ which gives:
$$\ker\Gamma=V\tilde{\uR}\EE_+(\tilde{\uR}V^*S V\tilde{\uR})=\EE_+(\uR^2S).$$
We have thus proved there is a smooth in $\zeta\neq0,$ matrix-valued
function $\uR,$ in the space of hermitian symmetric, positive
matrices such that:
$$\left[q\in\Ker \Gamma\right] \Leftrightarrow \left[\forall \zeta\neq0, \quad \uR^{-1}(\zeta)q\in \EE_+(\uR(\zeta) S(\zeta)\uR(\zeta))\right],$$
which is equivalent to: $$\forall \zeta\neq 0, \quad \Ker
\Gamma=\EE_+(\uR^2(\zeta) S(\zeta)).$$ Now consider $R$ defined, for
all $\zeta\neq 0,$ by:
$$R(\zeta)=\sqrt{S(\zeta)(\uR(\zeta))^2 S(\zeta)},$$
$$R(\zeta)=\sqrt{(\tilde{\uR}(\zeta)V^*(\zeta)S(\zeta))^*(\tilde{\uR}(\zeta)V^*(\zeta) S(\zeta))}.$$
$\zeta\mapsto R(\zeta)$ is smooth and, for all $\zeta,$ $R(\zeta)$
is a hermitian symmetric, positive definite matrix. Moreover, there
holds:
$$\left[q\in\Ker \Gamma\right] \Leftrightarrow \left[\forall \zeta\neq0, \quad R(\zeta)q\in \EE_+(R(\zeta) S^{-1}(\zeta)R(\zeta))\right],$$
which is equivalent to: $$\forall \zeta\neq0, \quad \Ker
\Gamma=\EE_+(S^{-1}(\zeta)R^2(\zeta)).$$ Let us detail the
computation of $R(\zeta).$
$$R(\zeta)=\sqrt{S(\zeta)V(\zeta)\tilde{\uR}^2(\zeta)V^*(\zeta) S(\zeta)}.$$
Moreover
$$\underline{\tilde{R}}^2(\zeta)=\tilde{S}^{-1}(\zeta) \tilde{R}^2
(\zeta)\tilde{S}^{-1}(\zeta),$$ we have thus:
$$R(\zeta)=\sqrt{\left(\tilde{R} (\zeta)\tilde{S}^{-1}(\zeta) V^*(\zeta) S(\zeta)\right)^*\left(\tilde{R} (\zeta)\tilde{S}^{-1}(\zeta) V^*(\zeta) S(\zeta)\right)},$$
which gives:
$$B(\zeta)=\left(\tilde{R} (\zeta)\tilde{S}^{-1}(\zeta) V^*(\zeta) S(\zeta)\right)^*\left(\tilde{R} (\zeta)\tilde{S}^{-1}(\zeta) V^*(\zeta) S(\zeta)\right).$$
We recall that $\tR$ is given, for all $\zeta\neq 0,$ by:
\begin{equation*}
\begin{aligned}{}
& \tR(\zeta) = \left[\begin{array}{cc}
Id_{N-p} & -\aleph(\zeta)\\
-\aleph^*(\zeta) & Id_{p}\\
\end{array}\right].
\end{aligned}
\end{equation*}
and that for all $\zeta\neq 0,$ $V(\zeta)$ is given by:
$$V(\zeta)= (\underline{\Lambda}(\zeta))^{-\frac{1}{2}} Z(\zeta),$$
where
$$\Lambda(\zeta)=Z^*(\zeta) \left(-S^{-1}\right)(\zeta)
Z(\zeta)$$ with $\Lambda$ is a diagonal matrix with real
coefficients: $(\lambda_1,\ldots,\lambda_N),$ and
$\underline{\Lambda}$ denotes the diagonal matrix with diagonal
coefficients $(|\lambda_1|,\ldots,|\lambda_N|).$
\begin{rem}
In the construction of $B$ the only freedom we have resides in the
choice of $\aleph.$
\end{rem}
\end{preuve}
\subsection{A change of dependent variables.}
\noindent Let us denote by $R:=B^{\frac{1}{2}}$ and $\hv:=R \hu.$
$\hv$ is hence solution of \eqref{deba}:
\begin{equation}
\left\{\begin{aligned}{} \label{deba}
& R^{-1}S R^{-1} \D_x \hv= R^{-1}SAR^{-1} \hv + R^{-1}S (A_d)^{-1}\hf, \quad \{x>0\},\\
& \Gamma R^{-1} \hv|_{x=0}=\Gamma \hat{g} ,\\
\end{aligned}\right.
\end{equation}
We will adopt the following notations: $S_R:=R^{-1}S R^{-1},$
$A_R:=RAR^{-1},$ and $\Gamma_R:=\Gamma R^{-1}.$ We first observe
that:
$$\ker \Gamma_R= R\ker\Gamma=R\EE_+((S)^{-1}R^2).$$
but $S_R^{-1}=RS^{-1}R$ thus
$$\ker \Gamma_R=R\EE_+(R^{-1}S_R R)=\EE_+(S_R).$$
This is where Lemma \ref{rau} is used in a crucial manner. Let us
denote by $\PP^-$ the projector on $\EE_-(S_R)$ parallel to
$\EE_+(S_R)$ and by by $\PP^+$ the projector on $\EE_+(S_R)$
parallel to $\EE_-(S_R);$ $\underline{\PP}^-$ and
$\underline{\PP}^+$ denoting the associated Fourier multiplier.
Since $S_R$ is hermitian symmetric, $\PP^-$ is in fact the
orthogonal projector on $\EE_-(S_R).$ The problem \eqref{deba} can
then be written:
\begin{equation*}
\left\{\begin{aligned}{}
& S_R \D_x \hv= S_R A_R \hv + R^{-1}S (A_d)^{-1}\hf, \quad \{x>0\},\\
& \PP^- \hv|_{x=0}=\PP^-\Gamma \hat{g} ,\\
\end{aligned}\right.
\end{equation*}
This problem is well-posed because, as a direct Corollary of
Proposition \ref{sym}, we have:
\begin{prop}\label{symr}
For all $\zeta$ such that $\tau^2+\gamma^2+ |\eta|^2=1,$ there
holds:
\begin{itemize}
\item $S_R(\zeta)$ is hermitian symmetric.
\item $\Re \left(S_R A_R\right)(\zeta)$ is positive definite.
\item $-S_R(\zeta)$ is definite negative on $\ker \Gamma_R$ and the dimension of $\ker
\Gamma_R$ is the same as the number of negative eigenvalues of
$-S_R(\zeta).$
\end{itemize}
\end{prop}

\begin{preuve}
For the sake of simplicity, let us omit the dependence in $\zeta$ in
our notations.
\begin{itemize}
\item $S_R:=R^{-1}S R^{-1},$ and both $S$ and $R$ are hermitian thus $S_R$ is hermitian.
\item $S_R A_R= R^{-1}SAR^{-1},$ thus for all $q\in\CC^N,$ there
holds: $$2\langle \Re(S_R A_R)q,q \rangle=\langle S_R A_Rq,q
\rangle+\langle q,S_R A_R q \rangle = \langle  R^{-1}SAR^{-1} q,q
\rangle+\langle q, R^{-1}SAR^{-1} q \rangle  ,$$ since $R^{-1}$ is
hermitian, we have then:
$$=\langle  SAR^{-1} q,R^{-1}q \rangle+\langle R^{-1}q, SAR^{-1} q \rangle= 2\langle \Re(SA)R^{-1}q,R^{-1}q\rangle.$$
Since $\Re(SA)$ is positive definite and $R$ is invertible, $\Re
\left(S_R A_R\right)$ is thus positive definite.

\item By construction of $R,$ it satisfies $\ker
\Gamma_R=\EE_+(S_R),$ with $S_R$ hermitian. As a consequence $-S_R$
is definite negative on $\ker \Gamma_R$ and the dimension of $\ker
\Gamma_R$ is the same as the number of negative eigenvalues of
$-S_R.$
\end{itemize}\end{preuve}
Let us mention that, since $R$ and $S$ remains uniformly bounded in
$\zeta\neq 0$, $\hf$ and $R^{-1}S (A_d)^{-1}\hf$ belongs to the same
space. In a same spirit as \cite{FG1}, this suggests the following
singular perturbation of \eqref{deba}:
\begin{equation*}
\begin{aligned}{}
&S_R \D_x \hv^\eps-\frac{1}{\eps} \PP^- \hv^\eps \textbf{1}_{x<0} = S_R A_R \hv^\eps - \frac{1}{\eps} \PP^- \Gamma \hat{g} \textbf{1}_{x<0}+ R^{-1}S (A_d)^{-1}\hf, \quad \{x\in\RR\},\\
\end{aligned}
\end{equation*}
This is equivalent to perturb \eqref{debm} as follow:
\begin{equation*}
\begin{aligned}{}
&S \D_x \hu^\eps-\frac{1}{\eps} R\PP^- R\hu^\eps \textbf{1}_{x<0} = S A \hu^\eps - \frac{1}{\eps} R\PP^-\Gamma \hat{g} \textbf{1}_{x<0}+ S (A_d)^{-1}\hf, \quad \{x\in\RR\},\\
\end{aligned}
\end{equation*}
Finally, this induces the following perturbation for \eqref{deb}:
\begin{equation}\label{pert}\left\{\begin{aligned}{}
&\cH \underline{u}^\eps+ \frac{1}{\eps} \MM \underline{u}^\eps \textbf{1}_{x<0} =f \textbf{1}_{x>0} + \frac{1}{\eps}\theta \textbf{1}_{x<0}, \quad \{x\in\RR\},\\
&\underline{u}^\eps|_{t<0}=0,\\
\end{aligned}\right.
\end{equation}
where $$\MM:= -e^{\gamma t} A_d \underline{S}^{-1} \underline{R} \underline{\PP}^-
\underline{R} e^{-\gamma t},$$ $$\theta=-e^{\gamma t} A_d
\underline{S}^{-1}\underline{R}\underline{\PP}^- \underline{\Gamma}
\tg,$$ and $\underline{S}(\D_t,\D_y)$ [resp
$\underline{R}(\D_t,\D_y)$] denotes the Fourier multiplier
associated to $S(\zeta)$ [resp $R(\zeta)$].

\section{Proof of Theorem \ref{kreiss}.}

First, we construct an approximate solution of equation \eqref{pert}
(which is also equation \eqref{penak}), then prove suitable energy
estimates that ensures $\underline{u}^\eps$ and its approximate
solution both converges towards the same limit as\newline
$\eps\rightarrow 0^+.$

\subsection{Construction of the approximate solution.}

$\uu^\eps$ is the solution of the well-posed Cauchy problem:
\begin{equation*}
\left\{\begin{aligned}{}
&\cH \underline{u}^\eps+ \frac{1}{\eps} \MM \underline{u}^\eps \textbf{1}_{x<0} =f\textbf{1}_{x>0}+ \frac{1}{\eps}\theta \textbf{1}_{x<0}, \quad \{x\in\RR\},\\
&\underline{u}^\eps|_{t<0}=0 .\\
\end{aligned}\right.
\end{equation*}
$\uu^\eps$ is moreover the solution of the well-posed Cauchy problem:
\begin{equation*}
\left\{\begin{aligned}{}
&\underline{S}A_d^{-1} \cH \underline{u}^\eps+ \frac{1}{\eps} \underline{S}A_d^{-1}\MM \underline{u}^\eps \textbf{1}_{x<0} =\underline{S}A_d^{-1}f \textbf{1}_{x>0}+\frac{1}{\eps} \underline{S}A_d^{-1}\theta \textbf{1}_{x<0}, \quad \{x\in\RR\},\\
&\underline{u}^\eps|_{t<0}=0 .\\
\end{aligned}\right.
\end{equation*}
The associated equation after tangential Fourier-Laplace transform
writes :
\begin{equation*}
\begin{aligned}{}
&S \D_x \uhu^\eps-\frac{1}{\eps} R\PP^- R\uhu^\eps \textbf{1}_{x<0} - S A \uhu^\eps = - \frac{1}{\eps} R\PP^-\Gamma \hat{g} \textbf{1}_{x<0}+ S (A_d)^{-1}\hf \textbf{1}_{x>0}, \quad \{x\in\RR\}.\\
\end{aligned}
\end{equation*}
or alternatively:
\begin{equation*}
\left\{\begin{aligned}{} &\uhu^\eps=R^{-1} \hv^\eps\\
&S_R \D_x \hv^\eps+\frac{1}{\eps} \PP^- \hv^\eps \textbf{1}_{x<0} = S_R A_R \hv^\eps + \frac{1}{\eps} \PP^-\Gamma \hat{g} \textbf{1}_{x<0}+ R^{-1}S (A_d)^{-1}\hf, \quad \{x\in\RR\},\\
\end{aligned}\right.
\end{equation*}
We will use the following formulation as a transmission problem in
our construction of an approximate solution:
\begin{equation*}
\left\{\begin{aligned}{}
&S_R \D_x \hv^{\eps+} = S_R A_R \hv^{\eps+} + R^{-1}S (A_d)^{-1}\hf, \quad \{x>0\},\\
&S_R \D_x \hv^{\eps-}+\frac{1}{\eps} \PP^- \hv^{\eps-} = S_R A_R \hv^{\eps-} + \frac{1}{\eps} \PP^-\Gamma \hat{g} , \quad \{x<0\},\\
&\hv^{\eps+}|_{x=0^+}=\hv^{\eps-}|_{x=0^-}.\\
\end{aligned}\right.
\end{equation*}

\noindent For $\Omega$ an open regular subset of $\RR^{d+1},$ and $\rho\in\NN,$ let us introduce the weighted spaces
$H_\gamma^{\varrho}(\Omega)$  defined by:
$$H_\gamma^{\varrho}(\Omega)=\{\varpi\in e^{\gamma t} L^2(\Omega), \|\varpi\|_{H_\gamma^{\varrho}(\Omega)}<\infty
\};$$ where
$$\|\varpi\|^2_{H_\gamma^{\varrho}(\Omega)}=\sum_{\alpha,|\alpha|\leq
\varrho} \gamma^{\rho-|\alpha|} \|e^{-\gamma t} \partial^\alpha
\varpi \|^2_ {L^2(\Omega)}.$$ We will construct an approximate
solution $\uu^\eps_{app}$ of $\uu^\eps.$ $\uu^\eps_{app}$ will be
constructed as follow:
$$\uu^\eps_{app}=\uu^{\eps+}_{app}\textbf{1}_{x>0}+\uu^{\eps-}_{app}\textbf{1}_{x<0},$$
where $\uu^{\eps\pm}_{app}$ is an approximate solution of
$\uu^{\eps\pm}$ satisfying the following ansatz:
$$\uu^{\eps\pm}_{app}=\sum_{j=0}^M \uU^\pm_j(\zeta,x) \eps^j,$$ where
the profiles $\uU^\pm_j$ belong to
$H^{k-\frac{3}{2}j}_\gamma(\Omega_T^\pm),$ where $\Omega_T^\pm$
stands for $[0,T]\times\RR^{d}_\pm.$ Denote
$$\hv^\eps_{app}=R \cF(e^{-\gamma t}
\uu^\eps_{app}):=\hv^{\eps+}_{app}\textbf{1}_{x>0}+\hv^{\eps-}_{app}\textbf{1}_{x<0}.$$
$\hv^{\eps\pm}_{app}$ is then an approximate solution of
$v^{\eps\pm}$ and is of the form:
$$\hv^{\eps\pm}_{app}=\sum_{j=0}^M V^\pm_j(\zeta,x) \eps^j;$$ where $$V^\pm_j=R \cF(e^{-\gamma t} \uU^\pm_j),$$ and conversely
$$ \uU^\pm_j=e^{\gamma t} \cF^{-1}\left(R^{-1}V^\pm_j\right).$$
The profiles $\uU^\pm_j$ can be constructed inductively at any order. Let us show how the first
profiles are constructed: Identifying the terms in $\eps^{-1}$
gives:
$$\PP^- V^-_0=\PP^-\Gamma \hat{g}.$$
Hence, $\PP^+ V^-_0$ remains to be computed in order to obtain the
profile $$ \uU^-_0=e^{\gamma t} \cF^{-1}\left(R^{-1}V^-_0\right).$$ Identifying the terms in $\eps^0$ gives then that
$V^+_0$ is solution of the well-posed problem:
\begin{equation}\label{ks1}
\left\{\begin{aligned}{}
&S_R \D_x V^+_0 = S_R A_R V^+_0 + R^{-1}S (A_d)^{-1}\hf, \quad \{x>0\},\\
&\PP^-V^+_0|_{x=0}=\PP^-\Gamma \hat{g}.\\
\end{aligned}\right.
\end{equation}
The associated profile $$ \uU^+_0=e^{\gamma t}
\cF^{-1}\left(R^{-1}V^+_0\right)$$ belongs then to
$H_\gamma^{k}(\Omega_T^+).$ Moreover, the problem \eqref{ks1} is
Kreiss-Symmetrizable and thus  the trace of the profile $\uU^+_0,$
see \cite{CP} for instance, satisfies:
$$\uU^+_0|_{x=0}\in H_\gamma^{k}(\Upsilon_T).$$
Since $V^+_0$ has just be computed, $V^-_0|_{x=0}$ is given by:  $V^+_0|_{x=0}-V^-_0|_{x=0}=0$ and thus, there holds:
$$\PP^-V^+_0|_{x=0}=\PP^-V^-_0|_{x=0}.$$ Moreover
$$S_R \D_x V^-_0-\PP^- V^-_1 = S_R A_R V^-_0 , \quad \{x<0\}.$$
Projecting this equation on $\EE_+(S_R)$ collinearly to $\EE_-(S_R)$
gives then:
$$S_R \D_x \PP^+ V^-_0- \PP^+ S_R A_R V^-_0 = 0, \quad \{x<0\},$$
Since $$\PP^+ S_R A_R V^-_0=\PP^+ S_R A_R \PP^+V^-_0+\PP^+ S_R A_R
\PP^-\Gamma \hat{g},$$ we have then:
$$S_R \D_x (\PP^+ V^-_0)- \PP^+ S_R A_R (\PP^+V^-_0)=\PP^+ S_R A_R
\PP^-\Gamma \hat{g}, \quad \{x<0\},$$
and as a consequence, $\PP^+ V^-_0$ is solution of the following
problem:
\begin{equation}\label{jus}
\left\{\begin{aligned}{} &S_R \D_x (\PP^+ V^-_0)- \PP^+ S_R A_R
(\PP^+V^-_0)=\PP^+ S_R A_R \PP^-\Gamma \hat{g} \quad \{x<0\},\\
&\PP^+V^-_0|_{x=0}=\PP^+V^+_0|_{x=0}.\\
\end{aligned}\right.
\end{equation}
Let us precise how \eqref{jus} has to
be interpreted: we denote $w=\PP^+ V^-_0.$ $w$ is then totally polarized
on $\EE_+(S_R),$ and satisfies the problem:
\begin{equation}\label{just}
\left\{\begin{aligned}{} &\PP^+ w=w\\
&S_R \D_x w- \PP^+ S_R A_R w=\PP^+ S_R A_R \PP^-\Gamma
\hat{g} \quad \{x<0\},\\
&w|_{x=0}=\PP^+V^+_0|_{x=0}.\\
\end{aligned}\right.
\end{equation}
As we will see, the problem $\eqref{just}$ is Kreiss-Symmetrizable
and thus well-posed. Indeed, for all $\zeta$ such that
$\tau^2+\gamma^2+ |\eta|^2=1,$ we have, omitting the dependencies in
$\zeta$ in our notations:
\begin{itemize}
\item For all $q\in\CC^N,$ there holds: $$\langle S_R
q,q \rangle=\langle q,S_R q \rangle.$$
\item Since $Re(S_R A_R)$ is positive definite and $\PP^+$ is an orthogonal projector, there is $C>0$ such
that, for all $q\in \EE_+(S_R),$ there holds: $$\langle \PP^+ S_R
A_R \PP^+ q,q \rangle+\langle q,\PP^+ S_R A_R \PP^+q \rangle \geq C
\langle q, q \rangle.$$ Indeed, for all $q\in \EE_+(S_R),$ there
holds:$$\langle \PP^+ S_R A_R \PP^+ q,q \rangle=\langle \PP^+ S_R
A_R \PP^+ q,\PP^+ q \rangle=\langle S_R A_R \PP^+ q,\PP^+ q
\rangle.$$
\item $-S_R$ is definite negative on $\ker \PP^+$ that is to say, that there is $c>0$ such that, for all $q\in \ker \PP^+,$ there holds: $$\langle -S_R q,q \rangle\leq -c\langle q, q \rangle.$$
Moreover $\ker \PP^+$ has the same dimension as the number of
negative eigenvalues in $-S_R.$
\end{itemize}
The profile $\uU^-_0$ can then be computed by:
$$ \uU^-_0:=e^{\gamma t} \cF^{-1}\left(R^{-1}(w+\PP^-\Gamma \hat{g})\right)$$ belongs to $H_\gamma^{k}(\Omega_T^-),$ moreover its trace $\uU^-_0|_{x=0}$
belongs \newline to $H_\gamma^{k}(\Upsilon_T).$ Consider now the
equation:
$$\PP^- V^-_1=S_R \D_x V^-_0 - S_R A_R V^-_0 , \quad \{x<0\}.$$
Since $\PP^- V^-_1|_{x=0}=\PP^- V^+_1|_{x=0},$ $V^+_1$ is solution of the well-posed problem:
\begin{equation*}
\left\{\begin{aligned}{}
&S_R \D_x V^+_1 = S_R A_R V^+_1, \quad \{x>0\},\\
&\PP^-V^+_1|_{x=0}=S_R \D_x V^-_0|_{x=0} - S_R A_R V^-_0|_{x=0}.\\
\end{aligned}\right.
\end{equation*}
Due to the loss of regularity in the boundary condition, the
associated profile $$ \uU^+_1=e^{\gamma t}
\cF^{-1}\left(R^{-1}V^+_1\right)$$ belongs to
$H_\gamma^{k-\frac{3}{2}}(\Omega_T^+),$ moreover its trace
$\uU^+_1|_{x=0}$ belongs \newline to
$H_\gamma^{k-\frac{3}{2}}(\Upsilon_T).$ Moreover, applying $\PP^+$
to the equation:
$$\PP^- V^-_2+ S_R A_R V^-_1=S_R \D_x V^-_1, \quad \{x<0\},$$
we obtain:
\begin{equation*}
\left\{\begin{aligned}{} &S_R \D_x \PP^+ V^-_1= \PP^+ S_R A_R
\PP^+V^-_1+\PP^+ S_R A_R \PP^-V^-_1 , \quad \{x<0\},\\
&\PP^+V^-_1|_{x=0}=\PP^+V^+_1|_{x=0}.\\
\end{aligned}\right.
\end{equation*}
As before, let us  take $\PP^+ V^-_1$  as the unknown of the well-posed problem:
\begin{equation*}
\left\{\begin{aligned}{}
&S_R \D_x (\PP^+ V^-_1)- \PP^+ S_R A_R
(\PP^+ V^-_1)=\PP^+ S_R A_R\left(S_R \D_x V^-_0 - S_R A_R V^-_0 \right), \, \, \{x<0\},\\
&(\PP^+ V^-_1)|_{x=0}=\PP^+V^+_1|_{x=0}.\\
\end{aligned}\right.
\end{equation*}
This problem is Kreiss-Symmetrizable since, for all $\zeta$ such
that \newline $\tau^2+\gamma^2+ |\eta|^2=1,$ there holds:
\begin{itemize}
\item For all $q\in\CC^N,$ there holds: $$\langle S_R
q,q \rangle=\langle q,S_R q \rangle.$$

\item There is $C>0$ such that
for all $q\in \EE_+(S_R),$ there holds: $$\langle \PP^+ S_R A_R
\PP^+ q,q \rangle+\langle q,\PP^+ S_R A_R \PP^+q \rangle \geq C
\langle q, q \rangle.$$

\item $-S_R$ is definite negative on $\ker \PP^+$ that is to say, that there is $c>0$ such that, for all $q\in \ker \PP^+,$ there holds: $$\langle -S_R q,q \rangle\leq -c\langle q, q \rangle.$$
Moreover $\ker \PP^+$ has the same dimension as the number of
negative eigenvalues in $-S_R.$

\end{itemize}
However, due to a loss of regularity in both the source term and the
boundary condition, the associated profile $$ \uU^-_1=e^{\gamma t}
\cF^{-1}\left(R^{-1}\left(\PP^+ V^-_1+S_R \D_x V^-_0 - S_R A_R V^-_0
\right)\right)$$ belongs to $H_\gamma^{k-\frac{3}{2}}(\Omega_T^-).$
The construction of the following profiles can be pursued at any
order the same way. In practice, we take:
$$u^{\eps+}_{app}=\uU_0^+,$$
$$u^{\eps-}_{app}=\uU_0^- +\eps \uU_1^-.$$
As a result, the approximate solution writes
$\uu^\eps_{app}:=\uu^{\eps+}_{app}\textbf{1}_{x>0}+\uu^{\eps-}_{app}\textbf{1}_{x<0};$
where $\uu^{\eps+}_{app}$ belongs to $H^{k}_\gamma(\Omega_T^+)$ and
$\uu^{\eps-}_{app}$ belongs to
$H^{k-\frac{3}{2}}_\gamma(\Omega_T^-).$ $\uu^\eps_{app}$ is then
solution of a well-posed problem of the form:
\begin{equation}\label{appp}
\left\{\begin{aligned}{}
&\cH \uu^{\eps}_{app}+ \frac{1}{\eps} \MM \uu^{\eps}_{app} \textbf{1}_{x<0} =f \textbf{1}_{x>0}+ \frac{1}{\eps}\theta \textbf{1}_{x<0}+\eps\ur^\eps, \quad \{x\in\RR\},\\
&\uu^{\eps}_{app}|_{t<0}=0 \quad .
\end{aligned}\right.
\end{equation}
Where
$\ur^\eps:=\ur^{\eps+}\textbf{1}_{x>0}+\ur^{\eps-}\textbf{1}_{x<0},$
with $\ur^{\eps+}\in H^{k-\frac{5}{2}}_\gamma(\Omega_T^+)$ and
\newline $\ur^{\eps-}\in H^{k-3}_\gamma(\Omega_T^-).$
\begin{rem}
In the case where $g=0,$ the loss of regularity in the profiles is
delayed by one step. As a result, in this case we obtain:
$$\uu^{\eps+}_{app}\in H^{k}_\gamma(\Omega_T^+),$$
$$\uu^{\eps-}_{app}\in H^{k}_\gamma(\Omega_T^-),$$
$$\ur^{\eps+}\in H^{k}_\gamma(\Omega_T^+),$$
$$\ur^{\eps-}\in H^{k-\frac{3}{2}}_\gamma(\Omega_T^-).$$
\end{rem}

\subsection{Stability estimates}
We will begin by proving energy estimates on the following equation:
\begin{equation}\label{simp}
\begin{aligned}{}
&S_R A_R \he^\eps -S_R \D_x \he^\eps+\frac{1}{\eps} \PP^- \he^\eps \textbf{1}_{x<0} =  \eps \hr^{\eps}, \quad \{x\in\RR\},\\
\end{aligned}
\end{equation}
where $\he^{\eps}=R\left(\cF(e^{-\gamma t}\uu^{\eps})-\cF(e^{-\gamma
t}\uu^{\eps}_{app})\right):=\hw^\eps;$ with
$w^\eps=\uu^{\eps}-\uu^{\eps}_{app}.$\newline Refering to
\eqref{appp}, $w^\eps$ is the solution of the Cauchy problem:
\begin{equation}\label{app1}
\left\{\begin{aligned}{}
& \cH w^{\eps}+\frac{1}{\eps} \MM w^{\eps} \textbf{1}_{x<0}= \eps \ur^\eps,\\
& w^{\eps}|_{t<0}=0 \quad .
\end{aligned}\right.
\end{equation}
\noindent For a fixed positive $\eps,$ the perturbation is
nonsingular and thus  the principal part of the pseudodifferential
operator $\cH+\frac{1}{\eps} \MM$ is the same as the principal part
of $\cH.$ Hence, there is a unique solution of the Cauchy problem
\eqref{app1}: $w^\eps$ which belongs to $H_\gamma^{k-3}(\Omega_T).$
In order to simplify the notations, in this chapter we shall denote
by $L^2$ and $H_\gamma^{\varrho}$ the spaces: $L^2(\Omega_T)$ and
$H_\gamma^{\varrho}(\Omega_T).$

\noindent We recall the definition of the weighted spaces:
$H_\gamma^{\varrho}(\Omega_T)$ for $\rho\in\NN.$
$$H_\gamma^{\varrho}(\Omega_T)=\{\varpi\in e^{\gamma t} L^2(\Omega_T), \|\varpi\|_{H_\gamma^{\varrho}(\Omega_T)}<\infty
\};$$ where
$$\|\varpi\|^2_{H_\gamma^{\varrho}(\Omega_T)}=\sum_{\alpha,|\alpha|\leq
\varrho} \gamma^{\rho-|\alpha|} \|e^{-\gamma t} \partial^\alpha
\varpi \|^2_ {L^2(\Omega_T)}.$$ For fixed positive $\eps,$ there
holds:
$$\int_{-\infty}^\infty \D_x\langle S_R \he^\eps,\he^\eps\rangle_{L^2} \, \, dx = 0.$$
$$\Leftrightarrow \int_{-\infty}^\infty 2Re\langle S_R \D_x\he^\eps,\he^\eps\rangle_{L^2} \, \, dx = 0.$$
Using the equation, we have then:
$$\int_{-\infty}^\infty Re\langle S_R A_R \he^\eps+\frac{1}{\eps} \PP^- \he^\eps - \eps \hr^{\eps},\he^\eps\rangle_{L^2} \, \, dx = 0.$$
which is equivalent to:
$$\int_{-\infty}^\infty Re\langle S_R A_R \he^\eps,\he^\eps\rangle_{L^2} \, \, dx= \frac{-1}{\eps}\int_{-\infty}^\infty Re\langle \PP^- \he^\eps  \eps \hr^{\eps},\he^\eps\rangle_{L^2} \, \, dx$$
$$+ \eps \int_{-\infty}^\infty Re\langle\hr^{\eps},\he^\eps\rangle_{L^2)} \, \, dx .$$
But $Re\langle S_R A_R \he^\eps,\he^\eps\rangle=\langle Re\left(S_R
A_R\right) \he^\eps,\he^\eps\rangle$ and $Re\left(S_R A_R\right)$ is
positive definite, hence there is $C>0,$ independent of $\eps$ such
that:
$$C\gamma\|\he^\eps\|^2_{L^2}+\frac{1}{\eps}\int_{-\infty}^{\infty} Re\langle \PP^-\he^\eps,\he^\eps\rangle \leq \int_{-\infty}^\infty Re\langle \eps\hr^{\eps},\he^\eps\rangle_{L^2} \, \, dx .$$
Thus, because $\PP^-$ is an orthogonal projector, for all positive
$\lambda,$ there holds:
$$C\gamma\|\he^\eps\|^2_{L^2}+\frac{1}{\eps} \|\PP^-\he^\eps\|^2_{L^2} \leq \frac{1}{2} \left( \frac{\gamma}{\lambda} \|\he^\eps\|^2_{L^2}+\frac{\lambda}{\gamma} \|\eps\hr^\eps\|^2_{L^2}\right).$$
Choosing $\lambda$ big enough we have $C-\frac{\eps}{2\lambda}>0$
and the following energy estimate:
$$\left(C-\frac{\eps}{2\lambda}\right)\gamma\|\he^\eps\|^2_{L^2}+\frac{1}{\eps} \|\PP^-\he^\eps\|^2_{L^2} \leq \frac{\eps^2\lambda}{2\gamma}  \|\hr^\eps\|^2_{L^2}.$$ This shows that $\he^\eps$ converges towards zero in $L^2$ when $\eps$ tends to zero, with a rate in $\cO(\eps).$
We recall that $\he^\eps$ is given by:
$$\he^\eps:=R\cF\left(e^{-\gamma t} (\uu^\eps_{app}-\uu^\eps)\right),$$
and $\hr^\eps$ is given by:
$$\hr^\eps:=R\cF\left(e^{-\gamma t} \ur^\eps \right).$$
Moreover, since $R$ and $\PP^-$ are two uniformly bounded, uniformly
definite positive matrices, there are two positive real numbers
$\alpha$ and $\beta$ such that, for all $\zeta\neq 0$ and $x\in\RR,$
there holds:
\begin{itemize}
\item $\alpha \|\cF\left(e^{-\gamma t}
(\uu^\eps_{app}-\uu^\eps)\right)\|^2_{L^2}  \leq
\|\he^\eps\|^2_{L^2}.$
\item $\alpha \|\PP^-\cF\left(e^{-\gamma t}
(\uu^\eps_{app}-\uu^\eps)\right)\|^2_{L^2}  \leq
\|\PP^-\he^\eps\|^2_{L^2}.$
\item $\|\hr^\eps\|^2_{L^2} \leq \beta
\|\cF\left(e^{-\gamma t} \ur^\eps\right)\|^2_{L^2}.$
\end{itemize}
Applying then Plancherel's equality we obtain then:
$$\left(C-\frac{\eps}{2\lambda}\right)\gamma\|\uu^\eps_{app}-\uu^\eps\|^2_{e^{\gamma t}L^2}+\frac{1}{\eps} \|\MM\left(\uu^\eps_{app}-\uu^\eps\right)\|^2_{e^{\gamma t}L^2} \leq \frac{\beta}{\alpha}\frac{\eps^2\lambda}{2\gamma}  \|\ur^\eps\|^2_{e^{\gamma t}L^2}.$$
We have thus proved there are two positive constants $c$ and $C$
such that:
$$c\gamma\|\uu^\eps_{app}-\uu^\eps\|^2_{e^{\gamma t}L^2}+\frac{1}{\eps} \|\MM\left(\uu^\eps_{app}-\uu^\eps\right)\|^2_{e^{\gamma t}L^2} \leq \frac{C\eps^2}{\gamma}  \|\ur^\eps\|^2_{e^{\gamma t}L^2}.$$
Let us denote by
$\|.\|^*_{H_\gamma^{\varrho}}:=\sqrt{\|.\|^2_{H_\gamma^{\varrho}(\Omega_T^+)}+\|.\|^2_{H_\gamma^{\varrho}(\Omega_T^-)}}.$
More generally, when $\ur^\eps\in H^{\varrho},$ there is two
positive constants $c_\rho$ and $C_\rho$ such that:
$$c_\rho \gamma\|\uu^\eps_{app}-\uu^\eps\|^{*2}_{H_\gamma^{\varrho}}+\frac{1}{\eps} \|\MM\left(\uu^\eps_{app}-\uu^\eps\right)\|^{*2}_{H_\gamma^{\varrho}} \leq \eps^2\frac{C_\rho}{\gamma}  \|\ur^\eps\|^{*2}_{H_\gamma^{\varrho}}.$$
As we have seen during the construction of the profiles,
$\varrho=k-3$ in general and $\rho=k-\frac{3}{2}$ in the case where
$g=0.$

\subsection{End of the proof of Theorem \ref{kreiss}.}
As a consequence of our stability estimate, there holds:
$$\|\uu^\eps_{app}-\uu^\eps\|^2_{H^{k-3}(\Omega_T^-)}+ \|\uu^\eps_{app}-\uu^\eps\|^2_{H^{k-3}(\Omega_T^+)}=\cO(\eps^2).$$
Moreover, by construction of $\uu^\eps_{app},$ there holds:
$$\|\uu^\eps_{app}-\underline{u}^-\|^2_{H^{k-3}(\Omega_T^-)}+ \|\uu^\eps_{app}-u\|^2_{H^{k-3}(\Omega_T^+)}=\cO(\eps^2).$$
Hence, we have proved that:
$$\|\uu^\eps-\underline{u}^-\|^2_{H^{k-3}(\Omega_T^-)}+ \|\uu^\eps-u\|^2_{H^{k-3}(\Omega_T^+)}=\cO(\eps^2).$$
By the same arguments, if $g=0,$ there holds:
$$\|\uu^\eps-\underline{u}^-\|^2_{H^{k-\frac{3}{2}}(\Omega_T^-)}+ \|\uu^\eps-u\|^2_{H^{k-\frac{3}{2}}(\Omega_T^+)}=\cO(\eps^2).$$
This completes the proof of Theorem \ref{kreiss}.

\section{Proof of Theorem \ref{maintheo}.}

Like in the proof of Theorem \ref{kreiss}, we begin by constructing
formally an approximate solution of equation \eqref{pena}. We prove
then suitable energy estimates that ensures both $u^\eps$ and its
approximate solution converges towards $\tu$ as $\eps\rightarrow
0^+.$

\subsection{Construction of an approximate solution.}

The goal of this Lemma is to replace the boundary condition $\Gamma
u|_{x=0}=\Gamma g$ of problem \eqref{deb} by a condition of the form
$\ubP^-(e^{-\gamma t}u)|_{x=0}=h$ with a suitable $h\in
H^k(\Upsilon_T).$
\begin{lem}\label{fund}
Let $u$ denote the unique solution in $H^k(\Omega_T^+)$ of the mixed
hyperbolic problem \eqref{deb}, $\ubP^+(\D_t,\D_y,\gamma)
\left(e^{-\gamma t}u\right)$ does not depend of the choice of the
boundary operator $\Gamma$ and of $g.$ Let us introduce the function
$h$ of $H^k(\Upsilon_T)$ defined by:
$$ \underline{\bP}^- \left(e^{-\gamma t}v|_{x=0}\right)+\underline{\bPi}\left(e^{-\gamma t}(g - v|_{x=0})\right).$$
The solution $u$ of the mixed hyperbolic problem \eqref{deb} is the
unique solution of the following well-posed mixed hyperbolic problem
\eqref{debr}:
\begin{equation}
\left\{\begin{aligned}{} \label{debr}
& \cH u=f, \quad \{x>0\},\\
& \ubP^-(\D_t,\D_y,\gamma) \left(e^{-\gamma t}u|_{x=0}\right)=h ,\\
& u|_{t<0}=0 \quad .
\end{aligned}\right.
\end{equation}
In addition, the mapping $(f,g)\rightarrow h$ is linear continuous
from \newline $H^k(\Omega_T^+)\times H^k(\Upsilon_T)$ to
$H^k(\Upsilon_T).$
\end{lem}
\begin{preuve}
Let $v$ denote a solution in $H^k(\Omega_T)$ of the equation:
\begin{equation*}
\left\{\begin{aligned}{}
& \cH v=f, \quad (t,y,x)\in\Omega_T,\\
& v|_{t<0}=0 \quad .
\end{aligned}\right.
\end{equation*}
We introduce then $\UU$ which is, by definition, the solution of the
following mixed hyperbolic problem:
\begin{equation*}
\left\{\begin{aligned}{}
& \cH \UU=0, \quad \{x>0\},\\
& \underline{\Gamma}(\D_t,\D_y,\gamma) \UU|_{x=0}=\underline{\Gamma}(\D_t,\D_y,\gamma) g - \underline{\Gamma}(\D_t,\D_y,\gamma) v|_{x=0},\\
& \UU|_{t<0}=0 \quad .
\end{aligned}\right.
\end{equation*}
The right hand side of the boundary condition is, a priori,
\newline in $H^{k-\frac{1}{2}}(\Upsilon_T).$ Hence the
solution $\UU$ belongs to $H^{k-\frac{1}{2}}(\Omega_T^+).$ By
construction we have:
\begin{equation}\label{eq}
\begin{aligned}{}
& u=\UU+v.\\
\end{aligned}
\end{equation}
Hence, since $u\in H^k(\Omega_T^+)$ and $v\in H^k(\Omega_T^+),$ in
fact we have: $$\UU\in H^k(\Omega_T^+).$$ Let $\hat{\UU}$ denote the
Fourier-Laplace transform in $(t,y)$ of $\UU$ (Fourier-Laplace
transform tangential to the boundary) given by: $\cF(e^{-\gamma t}
\UU).$ It satisfies the following symbolic equation:
\begin{equation*}
\left\{\begin{aligned}{}
& \D_x \hat{\UU}=A(\zeta)\hat{\UU}, \quad \{x>0\},\\
&\Gamma(\zeta) \hat{\UU}|_{x=0}=\Gamma(\zeta) \hat{g} - \Gamma(\zeta) \hv|_{x=0},\\
\end{aligned}\right.
\end{equation*}
where $\hat{g}$ and $\hv$ denotes respectively the tangential
Fourier-Laplace transform of $g$ and $v.$ Since $A(\zeta)$ is
independent of $x,$ projecting the above equation on
$\EE_+(A(\zeta))$ gives then:
$$\D_x \bP^+\hat{\UU}=A(\zeta)\bP^+\hat{\UU}.$$
Moreover $\bP^+\hat{\UU}|_{x=0}\in \EE_-(A(\zeta))\bigcap
\EE_+(A(\zeta))$ since $\lim_{x\rightarrow \infty}\bP^+\hat{\UU}=0.$
Hence, there holds:
$$\bP^+\hat{\UU}=0,$$ and thus $$\hat{\UU}=\bP^-\hat{\UU}.$$ The boundary condition:
$$\Gamma(\zeta) \hat{\UU}|_{x=0}=\Gamma(\zeta) \hat{g} -
\Gamma(\zeta) \hv|_{x=0}$$ is equivalent to:
$$\hat{\UU}|_{x=0}\in \hat{g} - \hv|_{x=0}+\Ker \Gamma.$$
We have thus:
$$\bP^-\hat{\UU}|_{x=0}\in \hat{g} - \hv|_{x=0}+\ker \Gamma.$$
Let us denote by $\bPi$ the projector on $\tilde{\EE}_-(A)$ parallel
to $\ker \Gamma,$ which has a sense because the Uniform Lopatinski
Condition holds. \newline Since $\hat{\UU}|_{x=0}\in
\tilde{\EE}_-(A),$ and of the Uniform Lopatinski Condition, the
above boundary condition is equivalent to:
$$\bPi \bP^-\hat{\UU}|_{x=0}= \bPi(\hat{g} - \hv|_{x=0}),$$
and thus, because $\bP^-\hat{\UU}|_{x=0}$ belongs to $\EE_-(A),$ we
have:
$$\bP^-\hat{\UU}|_{x=0}= \bPi(\hat{g} - \hv|_{x=0}).$$
As a consequence, we obtain:
$$\bP^-\hu|_{x=0}= \bP^-\hv|_{x=0}+\bPi(\hat{g} - \hv|_{x=0}).$$
Hence, there holds:
$$\underline{\bP}^- \left(e^{-\gamma t}u|_{x=0}\right)= \underline{\bP}^- \left(e^{-\gamma t}v|_{x=0}\right)+\underline{\bPi}\left(e^{-\gamma t}(g - v|_{x=0})\right):=h.$$
$\ubP^+(\D_t,\D_y,\gamma) \left(e^{-\gamma
t}u\right)=\ubP^+(\D_t,\D_y,\gamma) \left(e^{-\gamma t}v\right),$
thus it does not depend of the choice of the boundary operator
$\Gamma$ and of $g.$ Moreover, since $u|_{x=0}\in H^k(\Upsilon_T),$
it follows that \newline $g\in H^k(\Upsilon_T).$ Now, since the
Uniform Lopatinski Condition holds, $u$ satisfies the following
energy estimate:
$$\frac{1}{\gamma}\|u\|^2_{e^{\gamma t}L^2(\Omega_T^+)}+\|u|_{x=0}\|^2_{e^{\gamma t}L^2(\Upsilon_T)}\leq \gamma\|f\|^2_{e^{\gamma t}L^2(\Omega_T^+)}+\|g\|_{e^{\gamma t}L^2(\Upsilon_T)},$$
More generally, we have:
$$\frac{1}{\gamma}\|u\|^2_{H^k_\gamma(\Omega_T^+)}+\|u|_{x=0}\|^2_{H^k_\gamma(\Upsilon_T)}\leq \gamma\|f\|^2_{H^k_\gamma(\Omega_T^+)}+\|g\|_{H^k_\gamma(\Upsilon_T)}.$$
where $\|\varpi\|^2_{H^k_\gamma}:= \sum_{|\alpha|=0}^k
\gamma^{k-|\alpha|} \|\partial^\alpha \varpi\|^2_{e^{\gamma
t}L^2}.$\newline $h=\ubP^-(e^{-\gamma t}u|_{x=0})$ hence
$$\|h\|^2_{L^2(\Upsilon_T)}\leq C \|e^{-\gamma
t}u|_{x=0}\|^2_{L^2(\Upsilon_T)}= C \|u|_{x=0}\|^2_{e^{\gamma
t}L^2(\Upsilon_T)};$$ and for $0\leq j \leq d-1,$ there holds:
$$\|\D_j h\|^2_{L^2(\Upsilon_T)}\leq c_j\|\eta_j \bP^- \cF(e^{-\gamma t}u)|_{x=0} \|\leq
c'_j\|u|_{x=0}\|_{H^1_\gamma(\Upsilon_T)}.$$ More generally, we
have:
$$\|h\|^2_{H^k_\gamma(\Upsilon_T)}\leq  C_k\gamma\|f\|^2_{H^k_\gamma(\Omega_T^+)}+C_k\|g\|^2_{H^k_\gamma(\Upsilon_T)}.$$ But $\gamma$ is a positive real number fixed once and for all at the beginning of the paper, hence
this proves that the mapping $(f,g)\rightarrow h$ is continuous from
\newline $H^k(\Omega_T^+)\times H^k(\Upsilon_T)$
to $H^k(\Upsilon_T).$
\end{preuve}
As we will see, Lemma \ref{fund} is central in our construction of
an approximate solution. We will construct an approximate solution
$$u^{\eps}_{app}:=u^{\eps+}_{app}\textbf{1}_{x>0}+u^{\eps-}_{app}\textbf{1}_{x<0},$$
along the following ansatz:
$$u^{\eps+}_{app}:=\sum_{j=0}^M \eps^j u^+_j(t,y,x),$$ with $u^+_j\in H^{k-\frac{3}{2}j}_\gamma(\Omega_T^+),$ $u^+_j|_{x=0}\in H^{k-\frac{3}{2}j}_\gamma(\Upsilon_T);$ and
$$u^{\eps-}_{app}:=\sum_{j=0}^M \eps^j u^-_j(t,y,x),$$ with $u^-_j\in H^{k-\frac{3}{2}j}_\gamma(\Omega_T^-),$ $u^-_j|_{x=0}\in H^{k-\frac{3}{2}j}_\gamma(\Upsilon_T).$ As usual, we will refer to the terms $u^\pm_j$ as
profiles. We will rather work on the reformulation of
\newline problem \eqref{pena} as the transmission problem \eqref{tran}:
\begin{equation}\label{tran}
\left\{\begin{aligned}{}
& \cH u^{\eps+}=f, \quad \{x>0\},\\
& \cH u^{\eps-}+\frac{1}{\eps} A_d e^{\gamma t} \ubP^- e^{-\gamma t} u^{\eps-}=\frac{1}{\eps} A_d e^{\gamma t} \tilde{h}, \quad \{x<0\},\\
& u^{\eps+}|_{x=0}-u^{\eps-}|_{x=0}=0,\\
& u^{\eps\pm}|_{t<0}=0 \quad .
\end{aligned}\right.
\end{equation}
Plugging $u^{\eps+}_{app}$ and $u^{\eps-}_{app}$ in \eqref{tran} and
identifying the terms with same power in $\eps,$ we obtain the
following profiles equations:

\begin{itemize}
\item Identification of the terms of order $\eps^{-1}:$

\begin{equation} \label{-1-}
\begin{aligned}{} A_d e^{\gamma t} \ubP^- e^{-\gamma t}
u_0^-= A_d e^{\gamma t} \tilde{h}, \quad \{x<0\}. \end{aligned}
\end{equation}

\item Identification of the terms of order $\eps^{0}:$

\begin{equation} \label{0-}
\begin{aligned}{} \cH u_0^- + A_d e^{\gamma t} \ubP^- e^{-\gamma t}
u_1^-= 0, \quad \{x<0\}. \end{aligned}
\end{equation}

\begin{equation}\label{0+}
\begin{aligned}{}
& \cH u_0^{+}=f, \quad \{x>0\},\\
\end{aligned}
\end{equation}

\item Identification of the terms of order $\eps^{j}$ for $j\geq 1:$

\begin{equation} \label{j-}
\begin{aligned}{} \cH u_j^- + A_d e^{\gamma t} \ubP^- e^{-\gamma t}
u_{j+1}^-= 0, \quad \{x<0\}. \end{aligned}
\end{equation}

\begin{equation}\label{j+}
\begin{aligned}{}
& \cH u_j^{+}=0, \quad \{x>0\},\\
\end{aligned}
\end{equation}

\item Translation of the continuity condition over the boundary on the
profiles:

For all $1\leq j \leq M,$ there holds:
\begin{equation}\label{bc}
\begin{aligned}{}
& u_j^{+}|_{x=0}-u_j^{-}|_{x=0}=0.\\
\end{aligned}
\end{equation}
\end{itemize}
Denote by $\hu_{j}^\pm:=\cF(e^{-\gamma t} u_j^\pm)$ . We have then:
$$u_j^\pm:=e^{\gamma t} \cF^{-1}(\hu_{j}^\pm).$$ We will now
give the equations satisfied by the Fourier-Laplace transform of the
profiles: $\hu_{j}^\pm.$ First, equation \eqref{-1-} is equivalent
to:
\begin{equation*}
\begin{aligned}{} \bP^- \hu_0^-= \cF(\tilde{h}), \quad \{x<0\}. \end{aligned}
\end{equation*}
We deduce from this equation that there holds:
$$\bP^- \hu_0^-|_{x=0}= \hat{h}.$$
Then, using \eqref{bc} for $j=0,$ and \eqref{0+} gives that, for
$\gamma$ big enough, $$u_0^+=\cF(e^{-\gamma t} \hu_0^+),$$ where
$\hu_0^+$ is the solution of the well-posed first order ODE in $x$:
\begin{equation*}
\left\{\begin{aligned}{}
& \D_x \hu_0^{+}-A \hu_0^+= \cF(e^{-\gamma t}(A_d)^{-1}f), \quad \{x>0\},\\
& \bP^- \hu_0^+|_{x=0}= h,\\
\end{aligned}\right.
\end{equation*}
Thus $u_0^+$ is solution of:
\begin{equation*}
\left\{\begin{aligned}{}
& \cH u_0^{+}= f, \quad \{x>0\},\\
& e^{\gamma t}\ubP^- e^{-\gamma t}u_0^+|_{x=0}= h.\\
\end{aligned}\right.
\end{equation*}
Thanks to Lemma \ref{fund}, we recognize $u_0^+$ as the solution of
our starting mixed hyperbolic problem \eqref{deb}. Once $u_0^+$ is
known, so is $\hu_0^+$ and thus $\hu_0^-|_{x=0}$ is given by:
$$\hu_0^-|_{x=0}=\hu_0^+|_{x=0}.$$
Moreover, $$u_0^+|_{x=0}=u_0^-|_{x=0}\in H^k_\gamma(\Upsilon_T).$$
By \eqref{0-}, there holds:
\begin{equation*}
\begin{aligned}{} \D_x \hu_0^- -A \hu_0^- + \bP^- \hu_1^-= 0, \quad \{x<0\}. \end{aligned}
\end{equation*}
As a consequence, $\bP^+\hu_0^{-}$ is given by the well-posed ODE:
\begin{equation*}
\left\{\begin{aligned}{}
& \D_x (\bP^+\hu_0^{-})-A (\bP^+\hu_0^-)= 0, \quad \{x<0\},\\
& \bP^+\hu_0^-|_{x=0}= \bP^+\hu_0^+|_{x=0}.\\
\end{aligned}\right.
\end{equation*}
Indeed, since $\ker \bP^+(\zeta)=\EE_-(A(\zeta))$, this problem
satisfies the Uniform Lopatinski Condition: for all $\zeta\neq 0,$
there holds:
$$\EE_-(A(\zeta))\bigoplus\EE_+(A(\zeta))=\CC^N.$$
For $\gamma$ big enough, by linearity of the inverse Fourier
transform, $u_0^-$ can then be computed by:
$$u_0^-:=e^{\gamma t}\cF^{-1}(\bP^-\hu_0^-)+e^{\gamma t}\cF^{-1}(\bP^+\hu_0^-).$$
Following up with that process of construction, we can go on with
the construction of the profiles at any order. Indeed, assume that
all the profiles $(u_j^{+},u_j^{-})$ up to order $j$ have been
computed. Then consider the equation obtained through
identification:
\begin{equation*}
\begin{aligned}{} \bP^- \hu_{j+1}^-= -\D_x \hu_j^- +A \hu_j^- , \quad \{x<0\}. \end{aligned}
\end{equation*}
We see there is a loss of regularity between $\hu_{j+1}^-$ and
$\hu_j^-.$ \newline Let us say that $u_j^\pm\in
H^{m_j}_\gamma(\Omega_T^\pm).$ Considering the traces, we have then:
$u_j^\pm|_{x=0}\in H^{m_j-\frac{1}{2}}_\gamma(\Upsilon_T).$ We will
show in this part how the Sobolev regularity of the profiles
$u_{j+1}^\pm,$ which is by definition $m_{j+1},$ can be computed
knowing $m_j.$ To begin with $\ubP^- u_{j+1}^-$ belongs to
$H^{m_j-1}_\gamma(\Omega_T^-).$ $\ubP^- u_{j+1}^+|_{x=0},$ which
belongs to $H^{m_j-\frac{3}{2}}_\gamma(\Upsilon_T),$ is known by
$\ubP^- u_{j+1}^+|_{x=0}=e^{\gamma t} \cF^{-1}(\bP^-
\hu_{j+1}^+|_{x=0}),$ with:
$$\bP^-\hu_{j+1}^+|_{x=0}=\bP^- \hu_{j+1}^-|_{x=0}.$$ Hence, $\hu_{j+1}^+:=\cF(e^{-\gamma t}u_{j+1}^+)$
is the solution of the first order ODE in $x:$
\begin{equation*}
\left\{\begin{aligned}{}
& \D_x \hu_{j+1}^{+}-A \hu_{j+1}^+= 0, \quad \{x>0\},\\
& \bP^- \hu_{j+1}^+|_{x=0}= \bP^- \hu_{j+1}^-|_{x=0}.\\
\end{aligned}\right.
\end{equation*}
Since $\ker \bP^-(\zeta)=\EE_+(A(\zeta))$, this problem satisfies
the Uniform Lopatinski Condition: for all $\zeta\neq 0,$ there
holds:
$$\EE_-(A(\zeta))\bigoplus\EE_+(A(\zeta))=\CC^N.$$
As a consequence, this problem is well-posed and, $u_{j+1}^+\in
H_\gamma^{m_j-\frac{3}{2}}(\Omega_T^+).$ Moreover, there holds:
$$u_{j+1}^+|_{x=0}=u_{j+1}^-|_{x=0}\in H_\gamma^{m_j-\frac{3}{2}}(\Upsilon_T).$$
Indeed, $\bP^+ \hu_{j+1}^+\in H^{\infty}(\RR^{d+1}_+)$ hence $\ubP^+
u_{j+1}^+|_{x=0}\in H^{m_j-\frac{3}{2}}_\gamma(\Upsilon_T)$ and thus
$u_{j+1}^+|_{x=0}\in H^{m_j-\frac{3}{2}}_\gamma(\Upsilon_T).$
Furthermore, we have: $$u_{j+1}^-|_{x=0}=u_{j+1}^+|_{x=0}.$$
Applying $\bP^+$ on the following equation:
\begin{equation*}
\begin{aligned}{} \bP^- \hu_{j+2}^-= -\D_x \hu_{j+1}^- +A \hu_{j+1}^- , \quad \{x<0\}; \end{aligned}
\end{equation*}
we obtain then the equation:
\begin{equation*}
\begin{aligned}{} \D_x (\bP^+\hu_{j+1}^-) - A \bP^+\hu_{j+1}^-=0 , \quad \{x<0\}. \end{aligned}
\end{equation*}
\begin{rem}
\begin{equation*}
\begin{aligned}{} \bP^- \hu_{j+2}^-= -\D_x \hu_{j+1}^- +A \hu_{j+1}^- , \quad \{x<0\}. \end{aligned}
\end{equation*}
shows that the "Fourier profile" $\hu_{j+1}^-$ must be so that
$-\D_x \hu_{j+1}^- +A \hu_{j+1}^-$ is polarized on $\EE_-(A).$ It is
indeed the case because we search for $\hu_{j+1}^-$ satisfying:
\begin{equation*}
\begin{aligned}{} \D_x (\bP^+\hu_{j+1}^-) - A \bP^+\hu_{j+1}^-=0 , \quad \{x<0\}. \end{aligned}
\end{equation*}
\end{rem}
$u_{j+1}^-$ is given by:
$$u_{j+1}^-:=e^{\gamma t}\cF^{-1}(\bP^-\hu_{j+1}^-)+e^{\gamma t}\cF^{-1}(\bP^+\hu_{j+1}^-).$$
with $\ubP^+ u_{j+1}^-=e^{\gamma t}\cF^{-1}(\bP^+\hu_{j+1}^-)$
belongs to $H^{m_j-\frac{3}{2}}_\gamma(\Omega_T^+)$ and is the
unique solution of the well-posed first order ODE:
\begin{equation*}
\left\{\begin{aligned}{}
& \D_x (\bP^+\hu_{j+1}^{-})-A (\bP^+\hu_{j+1}^-)= 0, \quad \{x<0\},\\
& \bP^+ \hu_{j+1}^-|_{x=0}= \bP^+ \hu_{j+1}^+|_{x=0}.\\
\end{aligned}\right.
\end{equation*}
The profile $u_{j+1}^-$ belongs to
$H^{m_j-\frac{3}{2}}_\gamma(\Omega_T^-).$ This achieves to show that
the knowledge of $(u_{j}^+,u_{j}^-),$ allows us to compute
$(u_{j+1}^+,u_{j+1}^-).$ \newline Moreover
$m_{j+1}=m_j-\frac{3}{2},$ that is to say that a construction of
each supplementary profile consummate $\frac{3}{2}$ of Sobolev
regularity. In practice, we take:
$$u^{\eps+}_{app}=u_0^+,$$
$$u^{\eps-}_{app}=u_0^- +\eps u_1^-.$$
As a result, the approximate solution writes
$u^\eps_{app}:=u^{\eps+}_{app}\textbf{1}_{x>0}+u^{\eps-}_{app}\textbf{1}_{x<0};$
where $u^{\eps+}_{app}$ belongs to $H^{k}_\gamma(\Omega_T^+)$ and
$u^{\eps-}_{app}$ belongs to $H^{k-\frac{3}{2}}_\gamma(\Omega_T^-).$
The so defined $u^\eps_{app}$ is solution of a well-posed problem of
the form:
\begin{equation}\label{app}
\left\{\begin{aligned}{}
& \cH u^{\eps}_{app}+\frac{1}{\eps} A_d e^{\gamma t} \ubP^- e^{-\gamma t} u^{\eps}_{app} \textbf{1}_{x<0}= f\textbf{1}_{x>0}+\frac{1}{\eps} A_d e^{\gamma t} \tilde{h} \textbf{1}_{x<0}+ \eps r^\eps,\\
& u^{\eps}_{app}|_{t<0}=0 \quad .
\end{aligned}\right.
\end{equation}
Where $r^\eps:=r^{\eps+}\textbf{1}_{x>0}+r^{\eps-}\textbf{1}_{x<0},$
with $r^{\eps+}\in H^{k-\frac{5}{2}}_\gamma(\Omega_T^+)$ and
\newline $r^{\eps-}\in H^{k-3}_\gamma(\Omega_T^-).$

\subsection{Asymptotic Stability of the problem as $\eps$ tends towards zero.}
Denote by $v^\eps=u^{\eps}_{app}-u^\eps.$ By construction of
$u^{\eps}_{app},$ $v^\eps$ is solution of the following Cauchy
problem:
\begin{equation}\label{re}
\left\{\begin{aligned}{}
& \cH v^{\eps}+\frac{1}{\eps} A_d e^{\gamma t} \ubP^- e^{-\gamma t} v^{\eps} \textbf{1}_{x<0}=\eps r^\eps,\\
& v^{\eps}|_{t<0}=0 \quad .
\end{aligned}\right.
\end{equation}
For all positive $\eps,$ this problem is well-posed. In order to
investigate the stability of this problem as $\eps$ goes to zero, we
will reformulate it as a transmission problem. The restrictions of
$v^\eps$ to $\{x>0\}$ and $\{x<0\},$ respectively denoted by
$v^{\eps+}$ and $v^{\eps-}$ are solution the following transmission
problem:
\begin{equation}\label{transs}
\left\{\begin{aligned}{}
& \cH v^{\eps+}=\eps r^{\eps+}, \quad \{x>0\},\\
& \cH v^{\eps-}+\frac{1}{\eps} A_d e^{\gamma t} \ubP^- e^{-\gamma t} v^{\eps-}=\eps r^{\eps-}, \quad \{x<0\},\\
& v^{\eps+}|_{x=0}-v^{\eps-}|_{x=0}=0,\\
& v^{\eps\pm}|_{t<0}=0 \quad .
\end{aligned}\right.
\end{equation}
Let us denote by $V^\eps$ the function, valued in $\RR^{2N},$
defined for all $\{x>0\}$ and $(t,y)\in[0,T]\times \RR^{d-1}$ by:
$$V^\eps(t,y,x)=\left(\begin{array}{c}
V^{\eps+}(t,y,x)\\
V^{\eps-}(t,y,-x)\\
\end{array}\right) .$$
$v^\eps$ is solution of the Cauchy problem \eqref{re} iff $V^\eps$
is solution of the mixed hyperbolic problem on a half space
\eqref{eqstab} given below:
\begin{equation}
\left\{\begin{aligned}{}\label{eqstab}
& \tilde{\cH} V^{\eps}+B^\eps V^{\eps}=\eps R^\eps, \quad \{x>0\},\\
& \tilde{\Gamma} V^\eps|_{x=0}=0,\\
& V^\eps|_{t<0}=0 \quad ,
\end{aligned}\right.
\end{equation}
where $$\tilde{\cH}=\D_t+ \sum_{j=1}^{d-1} \left(\begin{array}{cc}
A_j & 0\\
0 & A_j\\
\end{array}\right)\D_j +\left(\begin{array}{cc}
A_d & 0\\
0 & -A_d\\
\end{array}\right)\D_x,$$
$$B^\eps=\left(\begin{array}{cc}
0 & 0\\
0 & \frac{1}{\eps} A_d e^{\gamma t} \ubP^- e^{-\gamma t}\\
\end{array}\right),$$
$$R^\eps(t,y,x)=\left(\begin{array}{c}
r^{\eps+}(t,y,x)\\
r^{\eps-}(t,y,-x)\\
\end{array}\right),$$
and
$$\tilde{\Gamma}=\left(\begin{array}{cc}
Id & -Id\\
\end{array}\right).$$

\noindent Returning to the construction of our approximate solution,
we have\newline $R^\eps\in
H^{k-\frac{5}{2}}_\gamma(\Omega_T^+)\times
H^{k-3}_\gamma(\Omega_T^+) $ and is such that $R^\eps|_{t<0}=0.$
\newline In fact $R^\eps\in H^{k'}_\gamma(\Omega_T^+)$ with $k'=k-3.$ For
all positive $\eps,$ there exists a unique solution $V^\eps$ in
$H^k_\gamma(\Omega_T^+)$ to the above problem. We will prove here
that this solution converges, uniformly in $\eps,$ towards $0$ in
$H^{k'}_\gamma(\Omega_T^+),$ as $\eps$ vanishes. As in the proof of
Kreiss Theorem, see \cite{CP} for instance, existence of solution
for mixed hyperbolic systems like \eqref{pena} or \eqref{eqstab},
are obtained through the proof of both direct and "dual" a priori
estimates on an adjoint problem. This estimates results in the
constant coefficient case of estimates on the Fourier-Laplace
transform of the solution. Additionally, if this "Fourier" estimate
can be proved, both direct and "dual" energy estimates are deduced
from it. In a first step, let us recall formally how to conduct the
Fourier-Laplace transform of a mixed hyperbolic problem:
\begin{equation*}
\left\{\begin{aligned}{}
&\cH u=f, \quad \{x>0\},\\
& \Gamma u|_{x=0}=g,\\
& u|_{t<0}=0 \quad ,
\end{aligned}\right.
\end{equation*}
Denote by $u_*:=e^{-\gamma t}u,$ $u_*$ is in particular a solution
of the following problem:
\begin{equation*}
\left\{\begin{aligned}{}
&\cH u_*+\gamma u_* =e^{-\gamma t}f, \quad \{x>0\},\\
&\Gamma u_*|_{x=0}=g \quad .
\end{aligned}\right.
\end{equation*}
We take then the tangential (with respect to (t,y)) Fourier
transform of this equation, which gives:
\begin{equation*}
\left\{\begin{aligned}{}
&A_d\D_x \hu_*+(\gamma+i\tau)\hu_*+ i\eta_j \sum_{j=1}^{d-1}A_j \hu_*=\cF\left(e^{-\gamma t}f\right), \quad \{x>0\},\\
&\Gamma \hu_*|_{x=0}=\hat{g} \quad .
\end{aligned}\right.
\end{equation*}
Multiplying this equation by $A_d^{-1},$ we obtain that $u^*$ is
solution of the following ODE in $x$:
\begin{equation*}
\left\{\begin{aligned}{}
&\D_x \hu_*-A \hu_*=(A_d)^{-1}\cF\left(e^{-\gamma t}f\right), \quad \{x>0\},\\
&\Gamma \hu_*|_{x=0}=\hat{g} \quad .
\end{aligned}\right.
\end{equation*}
Note that, $\hu_*$ and $u$ can be freely deduced from each other
through the formulas: $$\hu_*=\cF(e^{-\gamma t}u)$$ and
$$u=e^{\gamma t}\cF^{-1}(\hu_*).$$

We shall now introduce a rescaled solution $\uV^{\eps}$ of the
solution $V^{\eps}$ of \eqref{eqstab} defined as follows:
$\uV^{\eps}(t,y,x):=V^{\eps}(t,y,\eps x),$  and the rescaled
remainder: $\underline{R}^{\eps}(t,y,x):=R^{\eps}(t,y,\eps x).$
Denoting by $\underline{\hV}^{\eps}=\cF(e^{-\gamma
t}\underline{V}),$ the associated equation writes then:
\begin{equation*}
\left\{\begin{aligned}{}
& \D_x \underline{\hV}^{\eps}-\eps \tA \underline{\hV}^{\eps}+ M \underline{\hV}^{\eps}=\eps^2 \hat{R}^\eps, \quad \{x>0\},\\
& \tilde{\Gamma} \underline{\hV}^{\eps}|_{x=0}=0 \quad .
\end{aligned}\right.
\end{equation*}
where $$M(\zeta)=\left(\begin{array}{cc}
0 & 0\\
0 & \bP^-(\zeta)\\
\end{array}\right).$$
We remark that $$\eps \tA(\zeta)=\tA(\eps\zeta)= \tA(\hat{\zeta}),$$
with $\hat{\zeta}=(\hat{\tau},\hat{\gamma},\hat{\eta}) :=\eps\zeta.$
Moreover $\bP^-$ is homogeneous of order zero in $\zeta.$ Let us
denote $\tilde{R}^\eps(\hat{\zeta},x):=\hat{R}^\eps(\zeta,x)$ and
$\underline{\tV}^{\eps}(\hat{\zeta},x)
:=\underline{\hV}^{\eps}(\zeta,x).$ Hence $\underline{\tV}^{\eps}$
is solution of the following problem:
\begin{equation*}
\left\{\begin{aligned}{}
& \D_x \underline{\tV}^{\eps}+\left[-\tA(\hat{\zeta})+ M(\hat{\zeta})\right] \underline{\tV}^{\eps}=\eps^2 \tilde{R}^\eps(\hat{\zeta},x), \quad \{x>0\},\\
& \tilde{\Gamma} \underline{\tV}^{\eps}|_{x=0}=0 \quad .
\end{aligned}\right.
\end{equation*}
As a consequence, the Uniform Lopatinski Condition for problem
\eqref{eqstab} writes: For all $\hat{\gamma}>0,$
$$|det(\EE_-(\tA(\hat{\zeta})-M(\hat{\zeta}),\ker\Gamma )|\geq C>0.$$
In view of the proof of the Proposition \eqref{asympt}, we recall
that the spaces $\EE_\pm(A)$ have to be considered in the extended
sense defined above.
\begin{prop}\label{asympt}
Since $\cH$ satisfies the hyperbolicity Assumption in Assumption
\ref{vie}, the Uniform Lopatinski Condition is satisfied for our
present problem; that is to say that, for all $\hat{\zeta}$ such
that $\hat{\gamma}>0$ there holds:
$$|det(\EE_-(\tA(\hat{\zeta})-M(\hat{\zeta}),\ker\Gamma )|\geq C>0.$$
\end{prop}

\begin{preuve}
We will begin to show that the Uniform Lopatinski Condition writes
as well that for all $\hat{\zeta}\neq 0$ there holds:
\begin{equation}\label{lop}
\begin{aligned}{}
\EE_+(A(\hat{\zeta})-\bP^-(\hat{\zeta}))\bigcap
\EE_-(A(\hat{\zeta}))=\{0\} \quad .
\end{aligned}
\end{equation}
This notation keeps a sense for $\hat{\zeta}$ such that
$\hat{\gamma}=0$ because we will prove a posteriori that the
involved linear subspaces continuously extends from $\{\hat{\zeta},
\hat{\gamma}>0\}$ to $\{\hat{\zeta}, \hat{\gamma}=0\}.$  Then we
will prove that, for all $\hat{\zeta},$ the property \ref{lop} holds
true. The Uniform Lopatinski Condition writes actually, for all
$\hat{\zeta}\neq 0:$
$$\EE_-(\tA(\hat{\zeta})-M(\hat{\zeta}))\bigcap \ker\tilde{\Gamma}=\{0\}.$$
and thus, since we have:
$$\EE_-(\tA(\hat{\zeta})-M(\hat{\zeta}))=\EE_-(A(\hat{\zeta}))\times
\EE_+(A(\hat{\zeta})-\bP^-(\hat{\zeta})),$$ and by definition of
$\tilde{\Gamma},$ the Uniform Lopatinski Condition writes then that,
for all $\hat{\zeta}\neq 0,$ there holds:
$$\EE_+(A(\hat{\zeta})-\bP^-(\hat{\zeta}))\bigcap \EE_-(A(\hat{\zeta}))=\{0\}.$$
\begin{lem}\label{prel}
$$\EE_-\left(A(\hat{\zeta})-\bP^-(\hat{\zeta})\right)=\EE_-\left(A(\hat{\zeta})\right),$$
$$\EE_+\left(A(\hat{\zeta})-\bP^-(\hat{\zeta})\right)=\EE_+\left(A(\hat{\zeta})\right).$$
\end{lem}
\begin{preuve}
For all $\hat{\zeta}\neq 0,$ there is an invertible $N\times N$
matrix with complex coefficients $P(\hat{\zeta})$ such that:
$P^{-1}A P$ is trigonal and the diagonal coefficients are sorted by
increasing order of their real parts. Let us denote by $\nu$ the
dimension of $\EE_-\left(A\right).$ The above matrix $P$ traduces
the change of basis from the canonical basis of $\CC^N$ into
$(v_1,\ldots,v_\nu,v_{\nu+1},\ldots,v_N),$ where $$Span\left(
(v_k)_{1\leq k\leq\nu} \right)=\EE_-\left(A\right),$$ and
$$Span\left( (v_k)_{\nu+1 \leq k\leq N} \right)=\EE_+\left(A\right).$$
Moreover, there holds $$P^{-1}\bP^- P=D$$ where $D$ is the diagonal
matrix whose $\nu$ first diagonal terms are equal to $1$ and the
$N-\nu$ last diagonal terms are null.
$$P^{-1}(A-\bP^-) P=P^{-1}A P-D.$$ $P^{-1}A P-D$ is also trigonal, with the same eigenvalues with positive real part as $P^{-1}A P$ and the same eigenvalues with negative real part as $P^{-1}A P-Id.$
As a consequence, for all $\hat{\zeta}\neq 0,$ there holds:
$$\EE_-\left(A(\hat{\zeta})-\bP^-(\hat{\zeta})\right)=\EE_-\left(A(\hat{\zeta})\right),$$
$$\EE_+\left(A(\hat{\zeta})-\bP^-(\hat{\zeta})\right)=\EE_+\left(A(\hat{\zeta})\right).$$
\end{preuve}

As a consequence of Lemma \ref{prel}, the rescaled Uniform
Lopatinski Condition for $\eps>0, \eps \rightarrow 0$ happens to be
exactly the same as the one written for bigger positive $\eps.$
Indeed, it writes, for all $\hat{\zeta}\neq0:$
$$\EE_+(A(\hat{\zeta}))\bigcap \EE_-(A(\hat{\zeta}))=\{0\}.$$
\end{preuve}
The Lopatinski condition is satisfied, and, as a result, the
following, uniform in $\eps$, energy estimate holds for $\uV^\eps,$
for all $\gamma\geq \gamma_{k'}>0:$
\begin{equation*}
\begin{aligned}{} & \gamma \|\uV^\eps\|^2_{H_\gamma^{k'}(\Omega_T^+)}+\|\uV^\eps|^2_{x=0}\|_{H^{k'}_\gamma(\Upsilon_T)}\leq \frac{C}{\gamma}\|\eps
\underline{R}^\eps\|^2_{H^{k'}_\gamma(\Omega_T^+)} \quad ;
\end{aligned}
\end{equation*}
which is equivalent to:
\begin{equation}\label{estimate}
\begin{aligned}{} & \gamma \|V^\eps\|^2_{H^{k'}_\gamma(\Omega_T^+)}+\|V^\eps|_{x=0}\|^2_{H^{k'}_\gamma(\Upsilon_T)}\leq
\frac{C}{\gamma}\|\eps R^\eps\|^2_{H^{k'}_\gamma(\Omega_T^+)} \quad
.
\end{aligned}
\end{equation}
This proves the convergence of $V^\eps$ towards zero in
$H^{k'}_\gamma(\Omega_T^+).$ The weight $\gamma$ is fixed beforehand
thus, in fact, the solution of \eqref{eqstab} tends to zero in
$H^{k'}(\Omega_T^+)$ at a rate at least in $\cO(\eps).$

\section{End of proof of Theorem \ref{maintheo}.}
Let us consider $V^\eps$ defined by:
$$V^\eps(t,y,x):=\left(\begin{array}{c}
u^{\eps+}_{app}(t,y,x)-u^{\eps+}(t,y,x)\\
u^{\eps-}_{app}(t,y,-x)-u^{\eps-}(t,y,-x)\\
\end{array}\right).$$
This notation is perfectly fine because the so-defined function is
solution of an equation of the form \eqref{eqstab}. Moreover, thanks
to the stability estimate \eqref{estimate}, there is $\gamma_k$
positive such that, for all $\gamma>\gamma_k,$ we have:
$$\gamma\|u^{\eps}_{app}-u^{\eps}\|^2_{H^{k-3}_\gamma(\Omega_T^+)}+\gamma\|u^{\eps}_{app}-u^{\eps}\|^2_{H^{k-3}_\gamma(\Omega_T^-)}+\|u^{\eps}_{app}-u^{\eps}\|^2_{H^{k-3}_\gamma(\Upsilon_T)}\leq \frac{C}{\gamma}\|\eps R^{\eps+}\|^2_{H^{k-3}_\gamma(\Omega_T^+)}.$$
Hence, it follows that:
$$\|u^{\eps}_{app}-u^{\eps}\|^2_{H^{k-3}(\Omega^+_T)}+\|u^{\eps}_{app}-u^{\eps}\|^2_{H^{k-3}(\Omega^-_T)}=\cO(\eps^2).$$
Moreover, by construction of $u^{\eps}_{app},$ we have:
$$\|u^{\eps}_{app}-u\|^2_{H^{k-3}(\Omega^+_T)}+\|u^{\eps}_{app}-u^-\|^2_{H^{k-3}(\Omega^-_T)}=\cO(\eps^2).$$
As a result, we obtain that there holds:
$$\|u^{\eps}-u\|^2_{H^{k-3}(\Omega^+_T)}+\|u^{\eps}-u^-\|^2_{H^{k-3}(\Omega^-_T)}=\cO(\eps^2).$$
This concludes the proof of Theorem \ref{maintheo}.

\section{Appendix: answer to a question asked in \cite{paccou}.}
In this chapter, we will show that the loss of convergence observed
numerically in \cite{paccou} in a neighborhood of the boundary is
due to a boundary layer phenomenon. We consider the 1-D wave
equation:
\begin{equation}\label{wave}\left\{
\begin{aligned}{}
& \D_{tt}U-c^2\D_{xx} U=0, \quad (x,t)\in]0,\pi[\times\RR^+,\\
&U|_{x=0}=U|_{x=\pi}=0,\\
&U|_{t=0}(x)=sin(x),\\
&\D_tU|_{t=0}=0.\\
\end{aligned}\right.
\end{equation}
As in \cite{paccou}, we define then
$U^\eps=U^{\eps+}\textbf{1}_{x>0}+U^{\eps-}\textbf{1}_{x<0}$ by:
\begin{equation}\label{wavepen}\left\{
\begin{aligned}{}
& \D_{tt}U^{\eps+}-c^2\D_{xx} U^{\eps+}=0, \quad (x,t)\in]0,\pi[\times\RR^+,\\
& \D_{tt}U^{\eps-}-c^2\D_{xx} U^{\eps-}+ \frac{1}{\eps^2} U^{\eps-}=0, \quad (x,t)\in]-\infty,0[\times\RR^+,\\
&U^{\eps+}|_{x=0}-U^{\eps-}|_{x=0}=0\\
&\D_x U^{\eps+}|_{x=0}- \D_x U^{\eps-}|_{x=0}=0\\
&U^{\eps+}|_{x=\pi}=0.\\
&U^{\eps\pm}|_{t=0}(x)=sin(x), \quad\{\pm x>0\}.\\
&\D_t U^{\eps\pm}|_{t=0}=0, \quad\{\pm x>0\}.\\
\end{aligned}\right.
\end{equation}
We will now construct formally an approximate solution
$U^{\eps\pm}_{app}$ of $U^{\eps\pm}$ satisfying the following
ansatz:
$$U^{\eps+}_{app}=\sum_{j=0}^M U^+_j(t,x) \eps^j,$$
$$U^{\eps-}_{app}=\sum_{j=0}^M U^-_j\left(t,x,\left(\frac{x}{\eps}\right)\right) \eps^j,$$ where
the profiles $U^-_j(t,x,z):=\uU^-_j(t,x)+U^{*-}_j(t,z),$ with
$$\lim_{z\rightarrow -\infty} e^{-\alpha z}U^{*-}_j=0,$$ for some
$\alpha>0.$ Since the stability estimates are trivial here, we will
only focus on the construction of
$$U^{\eps}_{app}:=U^{\eps+}_{app}\textbf{1}_{x>0}+U^{\eps-}_{app}\textbf{1}_{x<0}.$$
Plugging $U^{\eps\pm}_{app}$ into problem \eqref{wavepen} and
identifying the terms with same power of $\eps,$ we obtain the
following equations: $$\uU^-_0=0,$$ Moreover, $U^{*-}_0=0$ as it is
the only solution of the problem:
\begin{equation*}\left\{
\begin{aligned}{}
& U^{*-}_0-c^2\D_{zz} U^{*-}_0=0, \quad \{z<0\},\\
& \D_z U^{*-}_0|_{z=0}=0,\\
& \lim_{z\rightarrow -\infty} U^{*-}_0=0.\\
\end{aligned}\right.
\end{equation*}
$U^{\eps+}_{app}$ converges towards $U^+_0$ as $\eps\rightarrow
0^+.$ As awaited $U^+_0$ is the solution of the well-posed 1-D wave
equation:
\begin{equation*}\left\{
\begin{aligned}{}
& \D_{tt}U^+_0-c^2\D_{xx} U^+_0=0, \quad (x,t)\in]0,\pi[\times\RR^+,\\
&U^+_0|_{x=0}=\uU^-_0|_{x=0}+U^{*-}_0|_{z=0}=0.\\
&U^+_0|_{x=\pi}=0.\\
&U^+_0|_{t=0}(x)=sin(x), \quad\{x>0\}.\\
&\D_t U^+_0|_{t=0}=0, \quad\{x>0\}.\\
\end{aligned}\right.
\end{equation*}
Let us write the following profiles equations: First, we can see
that, for all $j\geq 1,$ there holds:
$$\uU_j^-=0.$$
where $U^{*-}_1$ is the solution of the well-posed profile equation:
\begin{equation*}\left\{
\begin{aligned}{}
&  U^{*-}_1-c^2\D_{zz} U^{*-}_1=-\D_{tt} U^{*-}_{0}=0, \quad \{z<0\},\\
& \D_z U^{*-}_1|_{z=0}=\D_x U^+_0|_{x=0},\\
& \lim_{z\rightarrow -\infty} U^{*-}_1=0.\\
\end{aligned}\right.
\end{equation*}
Hence $U^{*-}_{1}$ is given by:
$$U^{*-}_{1}=c\D_x U^+_0|_{x=0} e^{\frac{z}{c}}.$$
We will show now that the profiles can be computed as any order.
Assume that $U^{*-}_{j}$ has been computed, $U^+_j$ is solution of
the well-posed 1-D wave equation:
\begin{equation*}\left\{
\begin{aligned}{}
& \D_{tt}U^+_j-c^2\D_{xx} U^+_j=0, \quad (x,t)\in]0,\pi[\times\RR^+,\\
&U^+_j|_{x=0}=U^{*-}_j|_{z=0}.\\
&U^+_j|_{x=\pi}=0.\\
&U^+_j|_{t=0}(x)=0, \quad\{x>0\}.\\
&\D_t U^+_0|_{t=0}=0, \quad\{x>0\}.\\
\end{aligned}\right.
\end{equation*}
$U^{*-}_{j+1}$ is then solution of the well-posed profile equation:
\begin{equation*}\left\{
\begin{aligned}{}
&  U^{*-}_{j+1}-c^2\D_{zz} U^{*-}_{j+1}=-\D_{tt} U^{*-}_{j}, \quad \{z<0\},\\
& \D_z U^{*-}_{j+1}|_{z=0}=\D_x U^+_{j}|_{x=0},\\
& \lim_{z\rightarrow -\infty} U^{*-}_{j+1}=0.\\
\end{aligned}\right.
\end{equation*}
Let us answer the question asked in \cite{paccou}: $U^{\eps-}$ is
bound to present boundary layer behavior in $\{x=0^-\},$ indeed its
approximate solution is composed \textbf{exclusively} of boundary
layer profiles, which describes quick transitions at the boundary
using a fast scale in $\eps$. As a result of the loss in convergence
induced by the boundary layer, the following estimate holds:
$$\|U^\eps-U\|_{L^2(]-\infty,\pi[\times\RR^+)}=\cO(\eps^{\frac{1}{2}}).$$
In \cite{paccou}, their small parameter is $\mu=\eps^2,$ as a
result, adopting the same notations as them, our estimate writes:
$$\|U^\mu-U\|_{L^2(]-\infty,\pi[\times\RR^+)}=\cO(\mu^{\frac{1}{4}}),$$
which is in agreement with the estimates given in \cite{paccou}.
Like in the penalization approach proposed by Bardos and Rauch
\cite{BR} and underlined by Droniou in \cite{D}, the boundary layer
only forms on one side of the boundary. The approximation
$U^{\eps+}$ of $U,$ is computed by taking
$U^{\eps+}|_{x=0}=U^{\eps-}|_{x=0},$ thus, in numerical
applications, the boundary layer phenomenon also affects the rate of
convergence of $U^{\eps+}$ towards $U,$\newline as $\eps\rightarrow
0^+.$

\newpage

\bibliographystyle{alpha}

\end{document}